\documentclass[12pt,reqno]{amsart}
\usepackage{appendix}
\usepackage[utf8]{inputenc}
\usepackage[T1]{fontenc}
\usepackage{lmodern}
\usepackage{amsmath,amssymb,amsfonts,amsthm,mathtools,color}
\usepackage[a4paper,left=2.0cm,right=2.0cm,top=1.8cm,bottom=2.0cm]{geometry}
\usepackage{microtype}
\usepackage{hyperref}

\numberwithin{equation}{section}

\newtheorem{theorem}{Theorem}[section]
\newtheorem{proposition}[theorem]{Proposition}
\newtheorem{lemma}[theorem]{Lemma}
\newtheorem{corollary}[theorem]{Corollary}
\newtheorem{remark}[theorem]{Remark}
\newtheorem{definition}[theorem]{Definition}
 
 \newcommand{\dis}[1]{}
 
\title[Sharp-interface limit of Navier--Stokes/Allen--Cahn system]
{The Sharp Interface Limit of the Two-Dimensional\\ Navier--Stokes/Allen--Cahn System as a Brakke Flow}
\author[Y. Liu]{Yuning Liu}
\address{NYU Shanghai, 567 Yangsi W road, Pudong, Shanghai 200126, China, and NYU-ECNU Institute of
Mathematical Sciences at NYU Shanghai, 3663 Zhongshan Road North, Shanghai, 200062, China}
\email{yl67@nyu.edu}
\author[Y. Tonegawa]{Yoshihiro Tonegawa}
\address{Department of Mathematics, Institute of Science Tokyo, 2-12-1 Ookayama, Meguro-ku, Tokyo 152-8551, Japan}
\email{tonegawa@math.sci.isct.ac.jp}
\date{\today}

\date{}
\begin{document}

\begin{abstract}
We study the sharp-interface limit of the two-dimensional
incompressible Navier--Stokes/Allen--Cahn system with constant mobility and prove the existence of global
integral varifold solutions to the coupled Navier--Stokes/mean curvature
flow system. The limiting interface measure satisfies Brakke's inequality, and its normal velocity equals the sum of its mean curvature
and the normal projection of the fluid velocity. 

\end{abstract}

\maketitle
 \tableofcontents

\section{Introduction and main theorem}

\subsection{Setting of the problem}
Diffuse-interface models provide a convenient description of two-phase fluid
flows without explicitly tracking a moving sharp boundary.  The two phases are
represented by a scalar order parameter which varies rapidly but smoothly
across a thin transition layer, so interfacial effects such as surface tension
are incorporated into the bulk equations.  Two of the most common phase-field
descriptions are based on the Allen--Cahn and Cahn--Hilliard equations.  The
former models local relaxation of the order parameter, while the latter
incorporates mass-conserving diffusion.  Coupling these phase-field equations
with the incompressible Navier--Stokes equations leads to diffuse-interface
models for viscous two-phase flow in which the fluid transports the phase
field and the diffuse interface exerts a capillary force on the fluid.

In the present work we consider the Navier--Stokes/Allen--Cahn system (with constant mobility) on the two dimensional flat torus
$\mathbb T^2= \mathbb R^2/(2\mathbb Z)^2$.  
Let
$$
W(s)=(1-s^2)^2,
\qquad
\sigma=\int_{-1}^1\sqrt{2W(s)}\,ds.
$$
For each $\epsilon>0$, consider the following
system of PDE
on $\mathbb T^2\times(0,\infty)$:
\begin{equation}\label{NSAC}
\left\{
\begin{aligned}
\partial_t\mathbf v_\epsilon
+(\mathbf v_\epsilon\cdot\nabla)\mathbf v_\epsilon
-\Delta\mathbf v_\epsilon+\nabla \Pi_\epsilon
&=-\epsilon\,\sigma^{-1}
\operatorname{div}(\nabla u_\epsilon\otimes\nabla u_\epsilon),
\\
\operatorname{div}\mathbf v_\epsilon&=0,
\\
\partial_tu_\epsilon+\mathbf v_\epsilon\cdot\nabla u_\epsilon
&=\Delta u_\epsilon-\frac1{\epsilon^2}W'(u_\epsilon),
\end{aligned}
\right.
\end{equation}
supplemented by suitable initial datum for
${\bf v}_\epsilon$ and $u_\epsilon$ at $t=0$.
When
the interfacial thickness $\epsilon>0$ tends to zero, one expects the diffuse layer to
collapse onto a sharp interface separating the two fluids.  In this limiting
description, the interface is represented by a moving curve, its motion
is determined jointly by the fluid transport and interfacial curvature, and
the diffuse capillary stress converges to a surface-tension force concentrated
on the interface. Indeed, a formal argument (assuming good convergence as $\epsilon\rightarrow 0$) shows that ${\bf v}=\lim_{\epsilon\rightarrow0}{\bf v}_\epsilon$
should satisfy
 \begin{equation}\label{sharp-system}
\left\{
\begin{aligned}
\partial_t \mathbf v+(\mathbf v\cdot\nabla)\mathbf v-\Delta \mathbf v+\nabla \Pi
&=H_{\mu_t}\,\sigma^{-1}\mu_t,
\\
\operatorname{div}\mathbf v&=0,
\\
\mathcal V&=H_{\mu_t}+P^\perp\mathbf v.
\end{aligned}
\right.
\end{equation}
 Here  $\mu_t$ is the one-dimensional measure of the interface, $H_{\mu_t}$ is the mean\footnote{Even though the dimension of the interface is one so that the reference ``mean curvature'' may not be proper, we use this terminology to be consistent with the higher dimensional language.} curvature vector of the interface, $\mathcal V$ is the normal velocity vector,
$P$ is the orthogonal projection onto the
tangent space $T_x\mu_t$, and $P^\perp=I-P$.
Understanding this limit is important both for justifying diffuse-interface
models from the viewpoint of continuum mechanics and for developing weak
formulations which remain meaningful when the sharp interface develops
singularities or changes topology.

\subsection{The statement of main result}
In the present paper, we give a rigorous proof on the convergence of \eqref{NSAC} to \eqref{sharp-system} for 
general initial datum, and by doing so,
prove the global-in-time existence of a weak solution of \eqref{sharp-system}, whose existence was 
not known before even for a short-time in this generality (for smooth initial datum
and the short-time existence of the regular solution, see Maekawa \cite{Maekawa}). The motion law of the 
interface is characterized in the sense of Brakke \cite{Brakke}, and the result is 
optimal for the natural energy class consideration. More precisely,
we characterize the limit of \eqref{NSAC}
and obtain the following result. 

\begin{theorem}\label{main-theorem}
Let $E_0\subset\mathbb T^2$ be a set of finite perimeter, and 
$$
\mathbf v_0\in L^2(\mathbb T^2;\mathbb R^2),
\qquad
\operatorname{div}\mathbf v_0=0 \text{ in distribution},
\qquad
\int_{\mathbb T^2}\mathbf v_0\,dx=0.
$$
  Then there
exist $(\mathbf v,\Pi)$ with  regularities 
\begin{align}\label{regularity v}
\mathbf v\in L^\infty([0,\infty);L^2(\mathbb T^2;\mathbb R^2))
\cap L^2
([0,\infty);\dot{H}^1(\mathbb T^2;\mathbb R^2)),\\
\Pi \in L_{\text {loc }}^2 ([0, \infty) ; L^r (\mathbb{T}^2 ) ) \quad \text { for every } r\in (1,2),
\end{align}
a family of Radon measures $\{\mu_t\}_{t\geq 0}$, the generalized mean
curvature $H_{\mu_t}\in L^2(d\mu_t\,dt;\mathbb R^2)$ for a.e.\,t>0, 
and a family of
sets of finite perimeter $E(t)\subset\mathbb T^2$ such that the following
properties hold.

\noindent
\emph{(1)}
The fluid velocity $\mathbf v$  satisfies $\operatorname{div}\mathbf v=0$ in distribution,
$\int_{\mathbb T^2}\mathbf v\,dx=0,$  and 
\begin{equation}\label{weak-NS}
\begin{aligned}
&-\int_0^\infty\int_{\mathbb T^2}
\mathbf v\cdot\partial_t\eta
+(\mathbf v\otimes\mathbf v):\nabla\eta\,
dx\,dt
+\int_0^\infty\int_{\mathbb T^2}
\nabla\mathbf v:\nabla\eta- \Pi\,{\rm div}\,\eta\,dx\,dt
\\
&\qquad=
\int_{\mathbb T^2}\mathbf v_0\cdot\eta(\cdot,0)\,dx
+\int_0^\infty\int_{\mathbb T^2}
H_{\mu_t}\cdot\eta\,d(\sigma^{-1}\mu_t)\,dt
\end{aligned}
\end{equation}
for every  
$\eta\in C_c^1(\mathbb T^2\times[0,\infty);\mathbb R^2)$.

\noindent
\emph{(2)}
For a.e.\ $t>0$, $\sigma^{-1}\mu_t$ is an integral one-dimensional
varifold, and the family
$\{\sigma^{-1}\mu_t\}_{t\in[0,\infty)}$ satisfies the following: for a.e.  $0\leq t_1<t_2<\infty$ and non-negative $\zeta\in C_c^1(\mathbb T^2\times[0,\infty))$, we have 
\begin{equation}\label{brakke-inequ}
    \left.\int_{\mathbb T^2}\zeta\,d\mu_t\right|_{t=t_1}^{t_2}\leq \int_{t_1}^{t_2}\int_{\mathbb T^2}\partial_t\zeta+(\nabla\zeta-H_{\mu_t}\zeta)\cdot(H_{\mu_t}+P^\perp{\bf v})\,d\mu_t dt,
\end{equation}
where   $P(x,t)$ is the orthogonal projection onto
$T_x\mu_t$, and $P^\perp=I-P$.
Moreover, $\mu_0=\sigma \|\nabla\mathcal X_{E_0}\|$,
$E(0)=E_0$ and 
\begin{equation}\label{phase-varifold}
\sigma\|\nabla\mathcal X_{E(t)}\|
\le\mu_t
\qquad\text{for all \ $t> 0$}.
\end{equation}

\noindent
\emph{(3)}
The map $t\mapsto\mathcal X_{E(t)}$ belongs to
$
C_{\mathrm{loc}}^{0,1/2}
\big([0,\infty);L^1(\mathbb T^2)\big)
$
and, for every $0\leq t_1<t_2<\infty$ and every
$\phi\in C^1(\mathbb T^2\times[0,\infty))$,
we have
\begin{equation}\label{def:bv-flow}
 \left.
\int_{E(t)}\phi(x,t)\,dx
\right|_{t=t_1}^{t_2}
 =
\int_{t_1}^{t_2}\int_{E(t)}
\left(
\partial_t\phi+\mathbf v\cdot\nabla\phi
\right)dx\,dt
+
\int_{t_1}^{t_2}\int_{\partial^*E(t)}
\phi\,H_{\mu_t}\cdot\nu_{E(t)}
\,d\mathcal H^1\,dt.
\end{equation}
Here, $\nu_{E(t)}$ is the outer-pointing unit normal of $E(t)$ which exists $\mathcal H^1$-a.e.\,on the reduced boundary $\partial^* E(t)$.

\end{theorem}
Note that the given initial velocity ${\bf v}_0$ and phase $E_0$ have minimal and natural 
regularity assumptions with respect to the energy law (see \eqref{enlaw}). The result claims the  global-in-time 
existence of four quantities, the fluid velocity field ${\bf v}$, the pressure $\Pi$, the interface measure $\mu_t$ and the phase $E(t)$. 
The property (1) is the weak formulation of the
Navier-Stokes equation with surface tension
effect.
The property (2) is a weak formulation for 
$\mathcal V=H_{\mu_t}+P^\perp{\bf v}$
in the sense of Brakke (see \cite[Chapter 2]{ToneBrakke}), and it specifies the rate of change of measure in time. The vector field  ${\bf v}$
is well-defined with respect to the surface measure $\mu_t dt$ as 
  $L^2$ functions (see Lemma \ref{limit-fluid-density-trace}). Even though \eqref{brakke-inequ} is proved for a.e. $0\leq t_1<t_2<\infty$, using a standard monotone argument (see   the proof of \cite[Prop 5.1]{LiuTonegawa}) we can upgrade it for all $0\leq t_1<t_2<\infty$. However, this involves modifications of $\{\sigma^{-1}\mu_t\}_{t\in[0,\infty)}$ at a countable set of  time. Note that such a family of measures comes directly from Allen-Cahn energy (see \eqref{con-mu_t} below).
The characterization (3)
expresses the same motion law of interface, but in
terms of the change of area $E(t)$ enclosed by the interface. More precisely, \eqref{def:bv-flow} states that the distributional time-derivative of $\mathcal X_{E(t)}$ equals $\nu_{E(t)}\cdot({\bf v}(\cdot,t)+H_{\mu_t})\mathcal H^1\llcorner_{\partial^* E(t)}$. 

\subsection{Some related works}
The sharp-interface limit of coupled Allen--Cahn/fluid systems has been
studied by several rather different methods.  For the present problem, it is
useful to distinguish three approaches: the varifold approach, which is
designed to produce global weak evolutions through singularities, and two
strong-solution approaches based respectively on spectral estimates and
relative entropy.

The geometric-measure-theoretic approach originates from the convergence of
Allen--Cahn equation to generalized mean curvature flow.  Ilmanen
\cite{Ilmanen} proved convergence of the Allen--Cahn equation to Brakke
flow, while Sato \cite{Sato} gave a short prove in
dimensions two and three utilizing R\"oger--Sch\"atzle's result \cite{RogerSchaetzle}.
For an Allen--Cahn equation with a prescribed transport field,
Liu--Sato--Tonegawa \cite{LiuSatoTonegawaTransport} proved, in dimensions
two and three, global existence of a generalized mean curvature flow with transport
term by establishing a uniform density-ratio bound for the diffuse surface
measures.  This was extended to higher dimensions by
Takasao--Tonegawa \cite{TakasaoTonegawa}.  Related diffuse-interface
estimates with Sobolev advection fields were developed by
Tonegawa--Tsukamoto \cite{TonegawaTsukamoto}.
See also Krämer--Laux--Takasao \cite{KramerLauxTakasao} for related work on mean curvature flow
with rough weights.

Along this line, most closely related to the present work is the paper of
Liu--Sato--Tonegawa \cite{LiuSatoTonegawaFluid}, where a genuinely coupled
two-phase fluid problem was treated by combining a Galerkin approximation
for the fluid with an Allen--Cahn phase-field approximation.  They obtained
global generalized solutions in dimensions two and three, with the phase boundary
moving by the sum of its mean curvature and the normal component of the
fluid velocity and exerting surface tension on the fluid.  Their fluid,
however, is shear-thickening and satisfies a power-law constitutive relation
with exponent
$
p>\frac{n+2}{2}.
$
Consequently, the Newtonian case $p=2$ is not covered: in two dimensions ($n=2$) it
is precisely the critical endpoint, while in three dimensions it lies below
the above threshold.
Another closely related work is carried out by Poiatti \cite{Poiatti2026}, which established the global existence
of De Giorgi-type varifold solutions to the coupled
Navier--Stokes/mean curvature flow system in dimensions two and three
by a minimizing movements approach. His solutions satisfy a
distributional transport equation for the phase indicator and a
sharp energy-dissipation inequality, rather than the geometric
motion law as an equality. The associated varifolds are shown to
be rectifiable and to possess square-integrable generalized mean
curvature, but their integrality remains open, as explicitly
pointed out in \cite[Remark 3.5]{Poiatti2026}. In contrast, our approach
constructs the sharp-interface limit directly from the coupled
Navier--Stokes/Allen--Cahn system and yields an integral varifold
together with the
characterization of 
motion in the sense of Brakke.

Even though it is not a coupled Navier-Stokes/Allen-Cahn system,
as a closely related work, Liu--Tonegawa \cite{LiuTonegawa} considered in
dimension two a prescribed transport field in the same critical Navier--Stokes
energy class
$$
\mathbf v\in
L^\infty_{\rm loc}([0,\infty);L^2(\mathbb R^2;\mathbb R^2))
\cap
L^2_{\rm loc}([0,\infty);H^1(\mathbb R^2;\mathbb R^2))
$$
and proved the  global-in-time existence of a curvature flow with forcing.  In that result the vector field ${\bf v}$ is given independently of
the evolving interface.  The difficulty in the present coupled problem is
that $\mathbf v_\epsilon$ itself depends on $\epsilon$ and is driven by the
diffuse surface tension force.

There is also a substantial literature on convergence to a prescribed
smooth sharp-interface evolution.  One approach is based on matched
asymptotic expansions together with spectral estimates for the linearized
Allen--Cahn operator.  For coupled fluid systems this program was initiated
by Abels--Liu \cite{AbelsLiu}, who treated the two-dimensional
Stokes/Allen--Cahn system.  Subsequent developments for
Navier--Stokes/Allen--Cahn systems include
\cite{AbelsFei,AbelsFeiLiuMoser,AbelsFeiMoser,AbelsMumtaz,
JiangSuXieSpectrum}.  These works cover, among other cases, two-dimensional
Navier--Stokes flow with equal or different viscosities, vanishing
mobility and mass-conserving Allen--Cahn dynamics; the recent work
\cite{AbelsFeiLiuMoser} extends the spectral approach to three-dimensional
Navier--Stokes/Allen--Cahn systems with nonconstant viscosity.

A second strong-solution approach is based on relative entropy developed by Fischer--Laux--Simon \cite{FischerLauxSimon} for Allen--Cahn equation.  For the
coupled Navier--Stokes/Allen--Cahn system this method was first used by
Hensel--Liu \cite{HenselLiu}, who proved convergence for constant mobility
and equal viscosities in both two and three dimensions, with quantitative
rates for well-prepared initial data as long as a sufficiently regular
sharp-interface solution exists.  
Substantial  extensions include  the
vanishing-mobility regime studied by Abels--Fischer--Moser
\cite{AbelsFischerMoser}.  

The spectral and relative-entropy approaches provide strong convergence and
quantitative information, but they rely on the existence of a sufficiently
regular limiting interface.  The varifold approach has a different purpose:
it seeks a  global-in-time weak evolution which remains meaningful after geometric
singularities occur.  The present work follows this latter direction and
addresses the Newtonian two-dimensional endpoint left open by
\cite{LiuSatoTonegawaFluid}.


\subsection{Ideas of proof}
Here we briefly describe the main idea of the proof. 
For a smooth solution to \eqref{sharp-system} one has
\begin{equation}
\label{enlaw}
\frac{d}{dt}
\left(
\frac12\int_{\mathbb T^2}|\mathbf v|^2\,dx
+\sigma^{-1}\mu_t(\mathbb T^2)
\right)
=
-\int_{\mathbb T^2}|\nabla\mathbf v|^2\,dx
-\int_{\mathbb T^2}|H_{\mu_t}|^2\,d(\sigma^{-1}\mu_t).
\end{equation}
Thus we have natural energy bounds 
which give
$$
\mathbf v\in L^\infty([0,\infty);L^2(\mathbb T^2))
\cap L^2
([0,\infty);\dot{H}^1(\mathbb T^2)),
\qquad
H_{\mu_t}\in L^2(d\mu_t\,dt).
$$
We have the similar $\epsilon$-version of \eqref{enlaw} (see \eqref{diffuse-energy}) which gives 
the similar energy bounds. 
Define
\begin{equation}
\label{stationary}
v_\epsilon=-\epsilon\Delta u_\epsilon+\frac1\epsilon W'(u_\epsilon),
\end{equation}
which may be considered as the chemical potential
in the context of thermodynamics of two-phase fluid.
The first difficulty in the study of $\epsilon$-limit of \eqref{NSAC}  is that the natural energy estimate controls the
material derivative
$$
\partial_tu_\epsilon
+
\mathbf v_\epsilon\cdot\nabla u_\epsilon
=
-\frac1\epsilon v_\epsilon,
$$
rather than $\partial_tu_\epsilon$ itself.  We shall prove a critical
one-dimensional density estimate for the diffuse surface measure and combine
it with the Meyers--Ziemer inequality to obtain
$$
\int_0^T\int_{\mathbb T^2}
|\mathbf v_\epsilon|^2\,d\mu_\epsilon^t\,dt\le C(T).
$$
See Proposition \ref{diffuse-density} in the sequel. 
This separates the transport term from the material derivative and yields
the Allen--Cahn action bound required by the compactness theory of
Mugnai--R\"oger \cite{MugnaiRoger,MugnaiRogerForced}.
Up to this point, we can characterize the 
limit interface as an $L^2$-flow. 

A second difficulty is identifying the limiting normal velocity of interface.
The energy estimates give strong convergence of the fluid velocity in
$L^2_{t,x}$, but not strong convergence of its gradient and hence not
directly strong convergence of its traces on the moving interfaces.  Here
the two-dimensional structure is essential.  Since the fluid velocity is
divergence free and has zero mean, we write
\begin{equation}\label{liu-revision-equ7}
\mathbf v_\epsilon=\nabla^\perp\Psi_\epsilon,
\qquad
\nabla^\perp\Psi_\epsilon=(-\partial_2\Psi_\epsilon,\partial_1\Psi_\epsilon),
\end{equation}
where $\Psi_\epsilon$ is the mean-zero scalar stream function. 
Strong $L^2$ convergence of $\mathbf v_\epsilon$ becomes strong $H^1$
convergence of $\Psi_\epsilon$.  The identity
$$
P_\epsilon^\perp\mathbf v_\epsilon
=P_\epsilon^\perp\nabla^\perp\Psi_\epsilon=
(P_\epsilon\nabla\Psi_\epsilon)^\perp
$$
converts the normal component of the fluid velocity into a tangential
derivative of the stream function.  Combining this identity with the first-variation formula for
the diffuse varifolds allows us to transfer the derivative to the test function and to identify the limiting motion law in the sense of $L^2$-flow,
\begin{equation}\label{motion-law-added}
\mathcal V
=
H_{\mu_t}+P^\perp\mathbf v.
\end{equation}
See Proposition \ref{normal-transport} below.

A third difficulty is to prove Brakke's inequality. 
The weak convergence of $P^\perp_\epsilon {\bf v}_\epsilon$ to $P^\perp{\bf v}$ is good enough to
give the $L^2$-flow characterization, which leads to Theorem 
\ref{main-theorem} (3), but not sufficient to prove
the inequality in (2) due to the coupling term
$H_{\epsilon}\cdot {\bf v}_\epsilon$. Instead of proving this 
directly, we first prove the localized energy inequality of \eqref{enlaw} first (see Proposition \ref{loc-ene-ineq} and \eqref{localized-coupled-energy-interval}). Note that we do 
not need to deal with the coupling term when $\epsilon\rightarrow 0$ for this energy identity.
The weak convergence of $P^\perp_\epsilon {\bf v}_\epsilon$ to $P^\perp{\bf v}$ already gives the
equality for the Navier-Stokes equation \eqref{weak-NS} and we may justify the 
substitution of the test function $\zeta{\bf v}$
in \eqref{weak-NS} by a suitable mollification. This gives the localized identity for the kinetic energy (see \eqref{weNS}). The subtraction of \eqref{weNS} from \eqref{localized-coupled-energy-interval}
gives the desired Brakke's inequality in (2).



 
\section{Critical diffuse density estimate and diffuse trace estimate}

Recall \eqref{NSAC} for the Navier--Stokes/Allen--Cahn system.
Following \cite[(2.1)--(2.4)]{RogerSchaetzle}, set
\begin{equation}\label{RS-v}
v_\epsilon
=-\epsilon\Delta u_\epsilon+\frac1\epsilon W'(u_\epsilon),
\end{equation}
and, for each time $t$,
\begin{equation}\label{RS-measures}
\begin{aligned}
\mu_\epsilon^t
&=\left(
\frac\epsilon2|\nabla u_\epsilon|^2
+\frac1\epsilon W(u_\epsilon)
\right)dx,
\\
\xi_\epsilon^t
&=\left(
\frac\epsilon2|\nabla u_\epsilon|^2
-\frac1\epsilon W(u_\epsilon)
\right)dx,
\\
\alpha_\epsilon^t
&=\frac1\epsilon v_\epsilon^2\, dx.
\end{aligned}
\end{equation}
Thus
\begin{equation}\label{gradient-measure}
\epsilon|\nabla u_\epsilon|^2 dx
=\mu_\epsilon^t+\xi_\epsilon^t.
\end{equation}

The initial data are chosen so that divergence-free smooth $\mathbf v_{\epsilon,0}$ satisfies
$$
\begin{aligned}
&\mathbf v_{\epsilon,0}\to\mathbf v_0
\quad\text{strongly in }L^2(\mathbb T^2),\\
&|u_{\epsilon,0}|\le1,
\quad
u_{\epsilon,0}\to 2\mathcal X_{E_0}-1\text{ strongly  in }L^1(\mathbb T^2),
\end{aligned}
$$
 and
\begin{equation}\label{well-prepared}
\lim_{\epsilon\to0}\mu_\epsilon^0(\mathbb T^2)
=\sigma\|\nabla\mathcal X_{E_0}\|(\mathbb T^2).
\end{equation}
Given a set of finite perimeter $E_0$, we may easily form
such a family $u_{\epsilon,0}\in C^\infty(\mathbb T^2)$ by approximations.

\begin{proposition}\label{fixed-epsilon}
Fix $\epsilon>0$ and let $m\ge3$. Assume
$$
\mathbf v_{\epsilon,0}\in H^m(\mathbb T^2;\mathbb R^2),
\qquad
\operatorname{div}\mathbf v_{\epsilon,0}=0,
\qquad
\int_{\mathbb T^2}\mathbf v_{\epsilon,0}\,dx=0,
$$
and
$$
u_{\epsilon,0}\in H^{m+1}(\mathbb T^2),
\qquad
|u_{\epsilon,0}|\le1.
$$
Then system \eqref{NSAC} admits a unique global-in-time  classical solution 
$(\mathbf v_\epsilon,u_\epsilon,\Pi_\epsilon)$ satisfying, for every $T>0$,
$$
\begin{aligned}
\mathbf v_\epsilon
&\in C([0,T];H^m(\mathbb T^2;\mathbb R^2))
   \cap L^2((0,T);H^{m+1}(\mathbb T^2;\mathbb R^2)),\\
u_\epsilon
&\in C([0,T];H^{m+1}(\mathbb T^2))
   \cap L^2((0,T);H^{m+2}(\mathbb T^2)), \\
   \Pi_\epsilon &  \in C([0,T];H^{m}(\mathbb T^2))\cap L^2((0,T);H^{m+1}(\mathbb T^2)).
\end{aligned}
$$
In particular, if the initial data are smooth, then
$\mathbf v_\epsilon$, $u_\epsilon$, and $\Pi_\epsilon$ are smooth on
$\mathbb T^2\times[0,\infty)$.
Moreover,
\begin{equation}\label{MP1}
|u_\epsilon(x,t)|\le1
\qquad\text{for all $(x,t)\in\mathbb T^2\times[0,\infty)$},
\end{equation}
and for all $t\ge0$,
\begin{equation}\label{diffuse-energy}
\begin{aligned}
& \frac12\|\mathbf v_\epsilon(t)\|_{L^2}^2
+\sigma^{-1}\mu_\epsilon^t(\mathbb T^2)
+\int_0^t\int_{\mathbb T^2}|\nabla\mathbf v_\epsilon|^2\,dx\,ds
+\sigma^{-1}\int_0^t\alpha_\epsilon^s(\mathbb T^2)\,ds
 \\
 &=
\frac12\|\mathbf v_{\epsilon,0}\|_{L^2}^2
+\sigma^{-1}\mu_\epsilon^0(\mathbb T^2).
\end{aligned}
\end{equation}

\end{proposition}

\begin{proof}
Existence and uniqueness of a classical solution follow from the
modified  Galerkin construction of 
Lin--Liu \cite{LinLiu}. The regularity of $\Pi_\epsilon$ follows from 
\[
-\Delta \Pi_\epsilon
=\partial_i\partial_j
\bigl(v_{\epsilon,i}v_{\epsilon,j}
+\sigma^{-1}\partial_i u_\epsilon\,\partial_j u_\epsilon\bigr),
\]
and that the operator $(-\Delta)^{-1}\partial_i\partial_j$ preserves Sobolev regularity.

The maximum principle applied to the Allen--Cahn equation gives
\eqref{MP1}. Indeed, testing the equation by $(u_\epsilon-1)_+$ gives
$$
\begin{aligned}
\frac12\frac{d}{dt}
\|(u_\epsilon-1)_+\|_{L^2}^2
&+
\|\nabla (u_\epsilon-1)_+\|_{L^2}^2
\\
&+
\frac1{\epsilon^2}
\int_{\{u_\epsilon>1\}}
W'(u_\epsilon)(u_\epsilon-1)\,dx
=0.
\end{aligned}
$$
Since $W'(s)(s-1)\ge0$ for $s\ge1$, one obtains
$u_\epsilon\le1$. The lower bound follows in the same way by testing with
$(-1-u_\epsilon)_+$.
\end{proof}

For each $t\geq 0$, define
$$
D_\epsilon(t)=
\sup_{x\in\mathbb T^2,\ 0<r\leq 1}
\frac{\mu_\epsilon^t(B_r(x))}{r}.
$$
The following Lemma \ref{monotone2} is proved in  
\cite[Proposition 4.7]{RogerSchaetzle},
taking into account that $|u_\epsilon|\leq 1$ and the dimension is two.
\begin{lemma}\label{monotone2}
Let $0<\beta\leq 1/2$ and $0<\epsilon\leq s\leq r\leq 1$, then
there exists $C$ depending only on $\beta$
such that
\begin{equation}\label{RS-quoted}
s^{-1}\mu^t_\epsilon(B_s(x))
\le C\big(1+r^{-1}\mu^t_\epsilon(B_r(x))+\alpha^t_\epsilon(B_{3r^{1-\beta}}(x))\big).
\end{equation}
\end{lemma}
With this estimate, we prove:
\begin{proposition}[Critical diffuse density bound]\label{diffuse-density}
There exists $C$ independent of $\epsilon$ such that for all
sufficiently small $\epsilon\in (0,1)$,
\begin{equation}\label{diffuse-density-est}
D_\epsilon(t)
\le C\left(
1+\mu_\epsilon^t(\mathbb T^2)
+\sqrt{\mu_\epsilon^t(\mathbb T^2)\,\alpha_\epsilon^t(\mathbb T^2)}
\right).
\end{equation}
Consequently, for every $T>0$,
\begin{equation}\label{density-L2}
\int_0^T D_\epsilon(t)^2\,dt\le C(T),
\end{equation}
where $C(T)$ is independent of $\epsilon$.
\end{proposition}

\begin{proof}
We divide the proof into scales above and below $\epsilon$.

\medskip
\noindent
\emph{Step 1: scales $r\ge\epsilon$.}
Let $x_0\in\mathbb T^2$ and let $R$ satisfy
\begin{equation}\label{constraint-R}
16\epsilon\le R\le 1.
\end{equation}
  In a flat coordinate ball centered at $x_0$, put
$$
\tilde u_{\tilde\epsilon}^t(y)=u^t_\epsilon(x_0+Ry),
\qquad
\tilde\epsilon=\frac\epsilon R,
\qquad
\tilde v^t_{\tilde\epsilon}(y)=R v^t_\epsilon(x_0+Ry).
$$
Then \eqref{stationary} becomes 
$$
-\tilde\epsilon\Delta\tilde u^t_{\tilde\epsilon}
+\frac1{\tilde\epsilon}W'(\tilde u^t_{\tilde\epsilon})
=\tilde v^t_{\tilde\epsilon},
$$
and the measures in \eqref{RS-measures} scale as
$$
\tilde\mu^t_{\tilde\epsilon}(A)
=R^{-1}\mu^t_\epsilon(x_0+RA),
\qquad
\tilde\alpha^t_{\tilde\epsilon}(A)
=R\alpha^t_\epsilon(x_0+RA).
$$

We apply \eqref{RS-quoted} with
$\beta=1/2$ and $r=1/16$.  Since $16\epsilon\le R$,
$\tilde\epsilon\le1/16$, while
$
3r^{1-\beta}=\frac34.
$
Since
$
\tilde\alpha^t_{\tilde\epsilon}(B_{3/4})
\le R\alpha^t_\epsilon(\mathbb T^2),
$
\eqref{RS-quoted} with $\tilde\epsilon$ gives the inequality
$$
\frac{\tilde\mu^t_{\tilde\epsilon}(B_\rho)}{\rho}
\le
C\,\big(1+\frac{\tilde\mu^t_{\tilde\epsilon}(B_{1/16})}{1/16}
+R\alpha^t_\epsilon(\mathbb T^2)\big),
\qquad
\tilde\epsilon\le\rho\le\frac1{16}.
$$
  Scaling back gives, for $\epsilon\le s\le R/16$,
\begin{equation}\label{RS-rescaled}
\frac{\mu^t_\epsilon(B_s(x_0))}{s}
\le C\big(
1+\frac{\mu^t_\epsilon(\mathbb T^2)}R
+R\alpha^t_\epsilon(\mathbb T^2)
\big).
\end{equation}
For $1\geq s>R/16$, the elementary bound
$$
\frac{\mu_\epsilon^t(B_s(x_0))}{s}
\le\frac{16\mu_\epsilon^t(\mathbb T^2)}R
$$
gives the same estimate after changing $C$.

If $\alpha_\epsilon^t(\mathbb T^2)>0$, set
$$
R_*=
\left(
\frac{\mu_\epsilon^t(\mathbb T^2)}
{\alpha_\epsilon^t(\mathbb T^2)}
\right)^{1/2}.
$$
The restriction imposed by 
$R$ in \eqref{constraint-R} leads us to consider three separate cases.
\begin{itemize}
\item 
If $16\epsilon\le R_*\le 1$, choose $R=R_*$ in
\eqref{RS-rescaled}. 
\item If $R_*<16\epsilon$, then for $s\ge\epsilon$,
$$
\frac{\mu_\epsilon^t(B_s(x_0))}{s}
\le\frac{\mu_\epsilon^t(\mathbb T^2)}\epsilon
<16\sqrt{\mu_\epsilon^t(\mathbb T^2)\alpha_\epsilon^t(\mathbb T^2)}.
$$
\item
If $R_*>1$, choose $R=1$ in
\eqref{RS-rescaled}; then
$\alpha_\epsilon^t(\mathbb T^2)<\mu_\epsilon^t(\mathbb T^2)$.
\end{itemize}

The case $\alpha_\epsilon^t(\mathbb T^2)=0$ is treated in the same way with
$R=1$.  We obtain
\begin{equation}\label{mesoscopic-density}
\sup_{x\in\mathbb T^2,\ \epsilon\le r\leq 1}
\frac{\mu_\epsilon^t(B_r(x))}{r}
\le C\left(
1+\mu_\epsilon^t(\mathbb T^2)
+\sqrt{\mu_\epsilon^t(\mathbb T^2)\alpha_\epsilon^t(\mathbb T^2)}
\right).
\end{equation}

\medskip
\noindent
\emph{Step 2: scales $r<\epsilon$.}
Fix $x_0$ and put
$$
U(y)=u_\epsilon(x_0+\epsilon y),
\qquad
V(y)=\epsilon v_\epsilon(x_0+\epsilon y).
$$
It satisfies
$$
-\Delta U+W'(U)=V\qquad\text{in }B_2.
$$
Set
$$
m=\int_{B_2}\left(\frac12|\nabla U|^2+W(U)\right)dy
=\frac{\mu_\epsilon^t(B_{2\epsilon}(x_0))}\epsilon,
$$
and
$$
k=\int_{B_2}|V|^2\,dy
\le\epsilon\alpha_\epsilon^t(\mathbb T^2).
$$
Since $|U|\le1$, the interior $W^{2,2}$ estimate gives
$$
\|\nabla^2U\|_{L^2(B_1)}\le C(1+k^{1/2}),
\qquad
\|\nabla U\|_{L^2(B_1)}^2\le Cm.
$$
By the two-dimensional Gagliardo-Nirenberg inequality, we have
\begin{equation*}
\begin{split}
\||\nabla U|^2\|_{L^2(B_1)} & \le C\left(\|\nabla|\nabla U|^2\|_{L^1(B_1)}+\||\nabla U|^2\|_{L^1(B_1)}\right) \\
& \le C\left(\|\nabla^2 U\|_{L^2(B_1)}\|\nabla U\|_{L^2(B_1)}+\|\nabla U\|_{L^2(B_1)}^2\right)
\le C\left(1+m+\sqrt{mk}\right).
\end{split}
\end{equation*}
For $0<\rho<1$, H\"older's inequality and boundedness of $W$ on $[-1,1]$
give
$$
\frac1\rho\int_{B_\rho}
\left(\frac12|\nabla U|^2+W(U)\right)dy
\le C\left(1+m+\sqrt{mk}\right).
$$
Scaling back, we obtain
\begin{equation}\label{subepsilon-pre}
\sup_{0<r<\epsilon}
\frac{\mu_\epsilon^t(B_r(x_0))}{r}
\le C\left(1+m+\sqrt{mk}\right).
\end{equation}
Taking $r=2\epsilon$ in  \eqref{mesoscopic-density} indicates that  
$
m\leq C
\big(1+\mu_\epsilon^t(\mathbb T^2)
+\sqrt{\mu_\epsilon^t(\mathbb T^2)\alpha_\epsilon^t(\mathbb T^2)}\big)
$
and we have $mk\leq \mu_\epsilon^t(B_{2\epsilon}(x_0))\alpha_\epsilon^t(\mathbb T^2)\leq \mu_\epsilon^t(\mathbb T^2)\alpha_\epsilon^t(\mathbb T^2)$.
Thus, \eqref{mesoscopic-density} and 
\eqref{subepsilon-pre} prove
\eqref{diffuse-density-est} on all scales.

Finally, \eqref{diffuse-energy} gives
$$
\sup_{0<t<T}\mu_\epsilon^t(\mathbb T^2)\le C,
\qquad
\int_0^T\alpha_\epsilon^t(\mathbb T^2)\,dt\le C.
$$
Squaring \eqref{diffuse-density-est} and integrating in time proves
\eqref{density-L2}.
\end{proof}

\begin{lemma}[Meyers--Ziemer   \cite{MeyersZiemer}]\label{trace-lemma}
Let $\lambda$ be a positive Radon measure on $\mathbb T^2$ satisfying
$$
\sup_{x\in\mathbb T^2,\ 0<r\leq 1}\frac{\lambda(B_r(x))}{r}\le D.
$$
Then every $f\in H^1(\mathbb T^2)$ satisfies
\begin{equation}\label{trace-est}
\int_{\mathbb T^2}|f|^2\,d\lambda
\le CD
\left(\|f\|_{L^2}\|\nabla f\|_{L^2}+\|f\|_{L^2}^2\right).
\end{equation}
\end{lemma}

\begin{proof} We extend $\mathbb T^2$ to $\mathbb R^2$ periodically and let $\phi\in C^\infty(\mathbb R^2;[0,1])$ be a cut-off function such that
$\phi=1$ on $B_2$ ($\mathbb T^2\subset B_2$) and $\phi=0$ outside of $B_3$. Consider $\tilde \lambda=\phi\lambda$ as a Radon measure on $\mathbb R^2$. Then, for any $B_r(x)\subset\mathbb R^2$, $0<r<\infty$, we have $\tilde\lambda(B_r(x))/r\leq C D$. 
Apply the Meyers--Ziemer inequality \cite{MeyersZiemer} to $\phi|f|^2$ on $\mathbb R^2$ and we obtain
$$\int_{\mathbb T^2}|f|^2\,d\lambda\leq \int_{B_3} \phi|f|^2\,d\tilde\lambda\leq CD\int_{B_3}|\nabla (\phi|f|^2)|\,dx\leq CD\int_{\mathbb T^2}|f|^2+|\nabla f||f|\,dx.
$$
By Cauchy-Schwarz, 
one obtains \eqref{trace-est}.
\end{proof}

\begin{proposition} \label{diffuse-trace}
For every $T>0$, the following inequalities hold uniformly in $\epsilon$:
\begin{equation}\label{u-trace}
\int_0^T\int_{\mathbb T^2}|\mathbf v_\epsilon|^2\,d\mu_\epsilon^t\,dt
\le C(T),
\end{equation}
\begin{equation}\label{time-curvature-bound}
\int_0^T\int_{\mathbb T^2}\epsilon|\partial_tu_\epsilon|^2\,dx\,dt
+\int_0^T\alpha_\epsilon^t(\mathbb T^2)\,dt
\le C(T).
\end{equation}
 
\end{proposition}

\begin{proof}
Apply Lemma \ref{trace-lemma} to $\lambda=\mu_\epsilon^t$ and
$f=|\mathbf v_\epsilon|$, using
$|\nabla|\mathbf v_\epsilon||\le|\nabla \mathbf v_\epsilon|$.  Proposition
\ref{diffuse-density} gives a uniform $L^2(0,T)$ bound for $D_\epsilon$,
while \eqref{diffuse-energy} gives uniform $L^\infty(0,T)$ bounds for
$\mu_\epsilon^t(\mathbb T^2)$ and $\|\mathbf v_\epsilon(t)\|_{L^2}$.  Hence
$$
\int_0^T\int|\mathbf v_\epsilon|^2\,d\mu_\epsilon^t\,dt
\le C\int_0^T
D_\epsilon(t)
\left(
\|\mathbf v_\epsilon\|_2\|\nabla\mathbf v_\epsilon\|_2
+\|\mathbf v_\epsilon\|_2^2
\right)dt
\le C(T)
$$
by Cauchy--Schwarz. 
The Allen--Cahn equation can be written, using \eqref{RS-v}, as
\begin{equation}\label{AC-v}
\partial_tu_\epsilon+\mathbf v_\epsilon\cdot\nabla u_\epsilon
=-\frac1\epsilon v_\epsilon.
\end{equation}
Combining \eqref{u-trace} with   \eqref{AC-v} and the curvature dissipation in
\eqref{diffuse-energy} leads to \eqref{time-curvature-bound}.
\end{proof}

\section{Compactness of varifold generated by $\mu_\epsilon^t$}
Define the diffuse
varifold
\begin{equation}\label{def:diffuse-varifold}
V_\epsilon^t(f)
=
\int_{\mathbb T^2}
f\big(x,P_\epsilon(x,t)\big)\,d\mu_\epsilon^t(x),
\qquad
f\in C^0(\mathbb T^2\times G(2,1)),
\end{equation}
where 
\begin{equation}
\label{projection-eps}
P_\epsilon=I-\nu_\epsilon\otimes\nu_\epsilon,
\qquad
P_\epsilon^\perp=\nu_\epsilon\otimes\nu_\epsilon,
\end{equation} with  $\nu_\epsilon=\nabla u_\epsilon/|\nabla u_\epsilon|$ on
$\{|\nabla u_\epsilon|>0\}$ and fixing $\nu_{\epsilon}$ to be an arbitrary unit vector on $\{|\nabla u_\epsilon|=0\}$ (for example $(0,1)$).
The first variation of \eqref{def:diffuse-varifold}  is deduced (\cite[Proposition 4.10]{RogerSchaetzle}) as
\begin{equation}\label{RS-first-variation}
(\delta V_\epsilon^t)( X )
=-\int_{\mathbb T^2}v_\epsilon\nabla u_\epsilon\cdot X \,dx
+\int_{\mathbb T^2}\nu_\epsilon^T\nabla X \,\nu_\epsilon\,d\xi_\epsilon^t.
\end{equation}

For later use define, on $\{|\nabla u_\epsilon|>0\}$,
\begin{equation}\label{def:diffusemcf}
H_\epsilon
=\frac{v_\epsilon}{\epsilon|\nabla u_\epsilon|}\nu_\epsilon,
\qquad
\mathcal V_\epsilon
=-\frac{\partial_tu_\epsilon}{|\nabla u_\epsilon|}\nu_\epsilon,
\end{equation}
and set them equal to zero elsewhere.  Since
$d(\mu_\epsilon^t+\xi_\epsilon^t)
=\epsilon|\nabla u_\epsilon|^2dx$, we can write 
\eqref{RS-first-variation} as
\begin{equation}\label{RS-first-variation+1}
(\delta V_\epsilon^t)( X )=-\int H_\epsilon\cdot X \,d(\mu_\epsilon^t+\xi_\epsilon^t)
+\int_{\mathbb T^2}\nu_\epsilon^T\nabla X \,\nu_\epsilon\,d\xi_\epsilon^t.
\end{equation}
The convective Allen--Cahn equation gives the exact identity
\begin{equation}\label{exact-diffuse-velocity}
\mathcal V_\epsilon
=H_\epsilon+P_\epsilon^\perp\mathbf v_\epsilon
\qquad
(\mu_\epsilon^t+\xi_\epsilon^t)dt\text{-a.e.}
\end{equation}
We   set
$$
d\mu_\epsilon=d\mu_\epsilon^t\,dt,
\qquad
d\xi_\epsilon=d\xi_\epsilon^t\,dt,
$$
and
\begin{equation}\label{def: mu-tilde}
d\widetilde\mu_\epsilon
:=
d(\mu_\epsilon+\xi_\epsilon)
=
\epsilon|\nabla u_\epsilon|^2\,dx\,dt.
\end{equation}
 
 Notice that
\begin{equation}\label{Heps-L2}
\int_{\mathbb T^2\times(0,T)}
|H_\epsilon|^2\,d\widetilde\mu_\epsilon
\le 
\int_0^T\alpha_\epsilon^t(\mathbb T^2)\,dt,
\end{equation}
whereas
\begin{equation}\label{Veps-L2}
\int_{\mathbb T^2\times(0,T)}
|\mathcal V_\epsilon|^2\,d\widetilde\mu_\epsilon
\le
\int_0^T\int_{\mathbb T^2}
\epsilon|\partial_tu_\epsilon|^2\,dx\,dt.
\end{equation}

Before stating the compactness result, we recall the meaning of convergence
of measure--function pairs due to Hutchinson \cite{Hutchinson}. 
\begin{definition}\cite[Section 4]{Hutchinson}
Let $\lambda_k,\lambda$ be positive Radon
measures on $\mathbb T^2\times[0,T]$, let
$f_k\in L^2(\lambda_k;\mathbb R^2)$ and
$f\in L^2(\lambda;\mathbb R^2)$. We write
$$
(\lambda_k,f_k)\to(\lambda,f)
$$
if $\lambda_k\to\lambda$ as Radon measures and
$$
\lim_{k\to\infty}
\int_{\mathbb T^2\times[0,T]}
f_k\cdot\zeta\,d\lambda_k
=
\int_{\mathbb T^2\times[0,T]}
f\cdot\zeta\,d\lambda
$$
for every
$\zeta\in C^0(\mathbb T^2\times[0,T];\mathbb R^2)$.
\end{definition}
\begin{lemma}\label{thmB.3}
\cite[Theorem 4.4.2]{Hutchinson} 
Assume that 
$
\sup_k\int |f_k|^2\,d\lambda_k<\infty,
$
then  
after passage to a subsequence,
$$
(\lambda_k,f_k)\to(\lambda,f)
$$
for some $f\in L^2(\lambda;\mathbb R^2)$, and 
$$
\int |f|^2\,d\lambda
\le
\liminf_{k\to\infty}\int |f_k|^2\,d\lambda_k.
$$
\end{lemma}

We next recall the $L^2$-flow formulation \cite[Definition 3.1]{MugnaiRoger}.
\begin{definition}\label{def-L2-flow}
Let $\{\lambda_t\}_{t\in(0,T)}$ be a family of Radon
measures on $\mathbb T^2$ such that $d\lambda=d\lambda_t\,dt$ is a Radon
measure on $\mathbb T^2\times(0,T)$.  Assume that, for a.e.\ $t$, $\lambda_t$ is a
rectifiable varifold and the associated varifold has the generalized mean curvature
$H_{\lambda_t}\in L^2(\lambda_t;\mathbb R^2)$.  The family
$\{\lambda_t\}_{t\in(0,T)}$ is an $L^2$-flow if there exist
$\mathcal V\in L^2(\lambda;\mathbb R^2)$ and $C_T<\infty$ such that
$$
\mathcal V(x,t)\perp T_x\lambda_t
\qquad\text{for $\lambda$-a.e.\ $(x,t)$},
$$
and
\[
\left|
\int_0^T\int_{\mathbb T^2}
\left(
\partial_t\zeta+ \mathcal V \cdot \nabla\zeta
\right)d\lambda
\right|
\le C_T\|\zeta\|_{C^0}
\]
for every $\zeta\in C_c^1(\mathbb T^2\times(0,T))$.
We call $\mathcal V$ the generalized velocity of $\lambda$. 
\end{definition}

The following result gives compactness properties of the quantities in \eqref{RS-measures}:
\begin{proposition}\label{MR-compactness}
For every $T>0$, there exist a sequence
$\epsilon_k\downarrow0$, a measurable set
$E\subset\mathbb T^2\times[0,T]$, a family of Radon measures
$\{\mu_t\}_{t\in[0,T]}$, a generalized mean curvature
$$
H_{\mu_t}\in L^2(d\mu_tdt;\mathbb R^2),
$$
and a generalized velocity
$$
\mathcal V\in L^2(d\mu_tdt;\mathbb R^2)
$$
such that the following properties hold.

\begin{enumerate}
\item $E$ has finite perimeter in
$\mathbb T^2\times[0,T]$ 
and
\begin{equation}
\label{uk-conv-char1}
u_{\epsilon_k}\to2\mathcal X_E-1
\qquad\text{strongly in }
L^1(\mathbb T^2\times[0,T]).
\end{equation}

\item
We have
\begin{equation}\label{con-mu_t}
\mu_{\epsilon_k}^t\to\mu_t
\qquad\text{as Radon measures on }\mathbb T^2
\end{equation}
for every $t\in[0,T]$ and
$$
\mu_{\epsilon_k}\to\mu:=\mu_tdt
\qquad\text{as Radon measures on }\mathbb T^2\times[0,T].
$$
For a.e.\ $t\in[0,T]$,
$\sigma^{-1}\mu_t$ is an integral one-dimensional varifold and
\begin{equation}\label{phase-measure-domination}
\|\nabla\mathcal X_{E(t)}\|=\mathcal H^1\llcorner\partial^*E(t)
\le
\sigma^{-1}\mu_t.
\end{equation}
where $E(t):=\{x\in\mathbb T^2:(x,t)\in E\}.$
\item
The discrepancy vanishes in total variation:
\begin{equation}\label{discrepancy-vanish}
|\xi_{\epsilon_k}|
\to0
\qquad\text{as Radon measures on }
\mathbb T^2\times[0,T].
\end{equation}
Consequently,
\begin{equation}\label{tilde-mu-limit}
\widetilde\mu_{\epsilon_k}
=
(\mu_{\epsilon_k}+\xi_{\epsilon_k})
\to\mu\text{ as Radon measures on }\mathbb T^2\times[0,T].
\end{equation}
 
\item
The diffuse mean curvatures converge to the generalized mean curvature in the
sense of measure--function pairs:
\begin{equation}\label{H-pair-convergence}
(\widetilde\mu_{\epsilon_k},H_{\epsilon_k})
\to
(\mu,H_{\mu_t}).
\end{equation}
Equivalently, for every
$\zeta\in C^0(\mathbb T^2\times[0,T];\mathbb R^2)$,
$$
\lim_{k\to\infty}
\int_{\mathbb T^2\times[0,T]}
H_{\epsilon_k}\cdot\zeta\,
d\widetilde\mu_{\epsilon_k}
=
\int_{\mathbb T^2\times[0,T]}
H_{\mu_t}\cdot\zeta\,d\mu.
$$

\item
The diffuse normal velocities converge in the sense of
measure--function pairs:
\begin{equation}\label{V-pair-convergence}
(\widetilde\mu_{\epsilon_k},\mathcal V_{\epsilon_k})
\to
(\mu,\mathcal V).
\end{equation}
Equivalently, for every
$\zeta\in C^0(\mathbb T^2\times[0,T];\mathbb R^2)$,
$$
\lim_{k\to\infty}
\int_{\mathbb T^2\times[0,T]}
\mathcal V_{\epsilon_k}\cdot\zeta\,
d\widetilde\mu_{\epsilon_k}
=
\int_{\mathbb T^2\times[0,T]}
\mathcal V\cdot\zeta\,d\mu.
$$
Moreover,
\begin{equation}
\label{orthogonality1}
\mathcal V(x,t)\perp T_x\mu_t
\qquad\text{for $\mu$-a.e.\ $(x,t)$}.
\end{equation}

\item
The family
$\{\sigma^{-1}\mu_t\}_{t\in(0,T)}$
is an $L^2$-flow with generalized velocity $\mathcal V$,
and 
\begin{equation}\label{curvature-liminf}
\int_{\mathbb T^2\times[0,T]}
|H_{\mu_t}|^2\,d\mu
\le
\liminf_{k\to\infty}
\int_0^T\alpha_{\epsilon_k}^t(\mathbb T^2)\,dt,
\end{equation}
\begin{equation}\label{velocity-liminf}
\int_{\mathbb T^2\times[0,T]}
|\mathcal V|^2\,d\mu
\le
\liminf_{k\to\infty}
\int_{\mathbb T^2\times[0,T]}
\epsilon_k|\partial_tu_{\epsilon_k}|^2\,dx\,dt.
\end{equation}
\end{enumerate}
\end{proposition}

\begin{proof}
The argument follows
\cite[Propositions 4.1--4.2, Theorems 4.3--4.4,
Proposition 6.1 and Lemmas 6.2--6.3]{MugnaiRoger}.
We give the details here.

\medskip
\noindent
\emph{Step 1: compactness of the phase fields.}
By \eqref{diffuse-energy} and \eqref{time-curvature-bound},
for every $T>0$,
\begin{equation}\label{MR-uniform-bounds}
\begin{aligned}
\sup_k\left[\sup_{0\le t\le T}\mu_{\epsilon_k}^t(\mathbb T^2)+\int_0^T
\alpha_{\epsilon_k}^t(\mathbb T^2)\,dt+\int_0^T\int_{\mathbb T^2}
\epsilon_k|\partial_tu_{\epsilon_k}|^2\,dx\,dt\right]&<\infty.
\end{aligned}
\end{equation}
The functioin 
$
\Phi(s):=\int_{-1}^s\sqrt{2W(r)}\,dr
$
satisfies 
$$
|\nabla\Phi(u_{\epsilon_k})|
=
\sqrt{2W(u_{\epsilon_k})}\,
|\nabla u_{\epsilon_k}|
\le
\frac{\epsilon_k}{2}|\nabla u_{\epsilon_k}|^2
+
\frac1{\epsilon_k}W(u_{\epsilon_k}),
$$
and similarly
$$
|\partial_t\Phi(u_{\epsilon_k})|
\le
\frac{\epsilon_k}{2}|\partial_tu_{\epsilon_k}|^2
+
\frac1{\epsilon_k}W(u_{\epsilon_k}).
$$
Thus $\{\Phi(u_{\epsilon_k})\}$ is bounded in
$BV(\mathbb T^2\times(0,T))$. Since
$$
\int_0^T\int_{\mathbb T^2}
W(u_{\epsilon_k})\,dx\,dt
\le
\epsilon_k T
\sup_{0<t<T}\mu_{\epsilon_k}^t(\mathbb T^2)
\to0,
$$
and $|u_{\epsilon_k}|\le1$, every $L^1$ limit of
$u_{\epsilon_k}$ takes only the values $-1$ and $1$.
Hence, after extracting a subsequence,
$$
u_{\epsilon_k}\to u=2\mathcal X_E-1\quad \text{ strongly in }L^1(\mathbb T^2\times(0,T)).
$$
After passing to another subsequence, Fubini's theorem gives
$$
u_{\epsilon_k}(\cdot,t)\to
2\mathcal X_{E(t)}-1
\quad\text{in }L^1(\mathbb T^2)\text{ for a.e. } t.
$$

\medskip
\noindent
\emph{Step 2: construction of the time slices $\mu_t$.}
For $\psi\in C^1(\mathbb T^2)$, differentiation of the diffuse energy gives
$$
\begin{aligned}
\frac{d}{dt}\mu_\epsilon^t(\psi)
&=
\int_{\mathbb T^2}
\psi\,v_\epsilon\,\partial_tu_\epsilon\,dx
-
\int_{\mathbb T^2}
\epsilon\,\partial_tu_\epsilon\,
\nabla u_\epsilon\cdot\nabla\psi\,dx.
\end{aligned}
$$
Consequently,
$$
\begin{aligned}
\int_0^T
\left|
\frac{d}{dt}\mu_\epsilon^t(\psi)
\right|dt
&\le
\|\psi\|_{C^0}
\left(
\int_0^T\alpha_\epsilon^t(\mathbb T^2)\,dt
\right)^{1/2}
\left(
\int_0^T\int_{\mathbb T^2}
\epsilon|\partial_tu_\epsilon|^2\,dx\,dt
\right)^{1/2}
\\
&\quad+
\|\nabla\psi\|_{C^0}
\left(
\int_0^T\int_{\mathbb T^2}
\epsilon|\partial_tu_\epsilon|^2\,dx\,dt
\right)^{1/2}
\left(
\int_0^T\int_{\mathbb T^2}
\epsilon|\nabla u_\epsilon|^2\,dx\,dt
\right)^{1/2}
\\
&\le C(T)\|\psi\|_{C^1}.
\end{aligned}
$$
Therefore $t\mapsto\mu_\epsilon^t(\psi)$ is uniformly bounded in
$BV(0,T)$.

Choose a countable dense set
$\{\psi_j\}_{j=1}^\infty\subset C^1(\mathbb T^2)$.
Helly's theorem and a diagonal argument give a subsequence
$\epsilon_k\downarrow0$ such that
$\mu_{\epsilon_k}^t(\psi_j)$ converges for every $j$ and every
$t$ outside a countable set $S$.  Since $S$ is countable, by another diagonal argument, we may assume that $\mu_{\epsilon_k}^t(\psi_j)$ converges on $S$ for every $j$. 
The uniform mass bound then defines a
Radon measure $\mu_t$ for every $t$ and
$$
\mu_{\epsilon_k}^t\to\mu_t\quad \text{ as Radon measures}.
$$
   By approximation of a general
$\eta\in C_c^0(\mathbb T^2\times(0,T))$ by finite sums
$\sum_j a_j(t)\psi_j(x)$, one obtains
$$
\mu_{\epsilon_k}\to\mu_tdt=\mu\quad \text{ as Radon measures}.
$$

\medskip
\noindent
\emph{Step 3:  compactness at a.e.\ time.}
Set
$$
g_k(t):=\alpha_{\epsilon_k}^t(\mathbb T^2)\ge0.
$$
By \eqref{MR-uniform-bounds} and Fatou's lemma,
$$
\int_0^T\liminf_{k\to\infty}g_k(t)\,dt
\le
\liminf_{k\to\infty}\int_0^Tg_k(t)\,dt<\infty.
$$
Hence
\begin{equation}
\liminf_{k\to\infty}
\alpha_{\epsilon_k}^t(\mathbb T^2)<\infty \quad \text{ for a.e.\ $t$.}
\end{equation}
Fix a time $t$ for which this is true. 
By \cite[Theorems 1.1]{MR2376729},
\cite[Theorems 4.1 and 5.1]{RogerSchaetzle}, we have 
\begin{equation}
\sigma^{-1}\mu_t
\quad\text{is an integral one-dimensional varifold}
\end{equation}
and $\mu_t$ has generalized mean curvature $H_{\mu_t}$ satisfying
\begin{equation}\label{static-H-bound}
\int_{\mathbb T^2}
|H_{\mu_t}|^2\,d\mu_t
\le
\liminf_{k\to\infty}
\alpha_{\epsilon_k}^t(\mathbb T^2).
\end{equation}
The Modica--Mortola lower bound at the same time gives
$$
\frac{\sigma}{2}
\|\nabla(2\mathcal X_{E(t)}-1)\|
\le\mu_t.
$$
Since
$
\|\nabla(2\mathcal X_{E(t)}-1)\|
=
2\mathcal H^1\llcorner\partial^*E(t),
$
we obtain \eqref{phase-measure-domination}.

We now integrate \eqref{static-H-bound} and use 
Fatou's lemma:
$$
\begin{aligned}
\int_0^T\int_{\mathbb T^2}
|H_{\mu_t}|^2\,d\mu_t\,dt
&\le
\int_0^T
\liminf_{k\to\infty}
\alpha_{\epsilon_k}^t(\mathbb T^2)\,dt
\\
&\le
\liminf_{k\to\infty}
\int_0^T
\alpha_{\epsilon_k}^t(\mathbb T^2)\,dt.
\end{aligned}
$$
This proves \eqref{curvature-liminf}.

\medskip
\noindent
\emph{Step 4: vanishing of the discrepancy.}
For $L>0$ define the bad-time set
$$
\mathcal B_{\epsilon,L}
:=
\left\{
t\in(0,T):
\alpha_\epsilon^t(\mathbb T^2)>L
\right\}.
$$
By Chebyshev's inequality and \eqref{MR-uniform-bounds},
$
|\mathcal B_{\epsilon,L}|
\le\frac{C}{L}.
$
For fixed $L$, the static discrepancy theorem implies
\begin{equation}\label{truncated-discrepancy}
\mathbf 1_{(0,T)\setminus\mathcal B_{\epsilon_k,L}}(t)
\,|\xi_{\epsilon_k}^t|(\mathbb T^2)\longrightarrow0
\qquad\text{for  every  \ }t\in(0,T).
\end{equation}
Indeed, if this failed at a fixed time, there would be a subsequence
with diffuse curvature bounded by $L$ and discrepancy bounded away from
zero, contradicting the static theorem.  The uniform mass bound and
$|\xi_{\epsilon_k}^t|\le\mu_{\epsilon_k}^t$ give an integrable
dominating constant.  Thus the integral in time of the left-hand side of
\eqref{truncated-discrepancy} tends to zero.
On the bad times,
\begin{equation}\label{bad-time-discrepancy}
\int_{\mathcal B_{\epsilon_k,L}}
|\xi_{\epsilon_k}^t|(\mathbb T^2)\,dt
\le \sup_{k,t}\mu_{\epsilon_k}^t(\mathbb T^2)
|\mathcal B_{\epsilon_k,L}|\le\frac{C(T)}L.
\end{equation}
Taking first $k\to\infty$ and then $L\to\infty$ proves, in fact,
\begin{equation}\label{discrepancy-total-mass}
\int_0^T|\xi_{\epsilon_k}^t|(\mathbb T^2)\,dt\longrightarrow0.
\end{equation}
This proves \eqref{discrepancy-vanish}; \eqref{def: mu-tilde} also gives
\eqref{tilde-mu-limit}.

\medskip
\noindent
\emph{Step 5: space-time varifold convergence and a common subsequence.}
This step uses only Steps 2--3 and the static varifold compactness
theorem, and does not use the curvature identification in Step 6.
For a.e. $t$, let $V^t$ be the rectifiable varifold with weight $\mu_t$
and tangent projection $P(x,t)$.  We use measurable representatives of
this family and of $P$; these can be chosen since $t\mapsto\mu_t$ is
measurable and the approximate spatial tangent is determined by $\mu_t$.
Set
\begin{equation}\label{varifold-uniform-mass}
M:=\max\left\{1,\sup_{k,t}\mu_{\epsilon_k}^t(\mathbb T^2)\right\},
\qquad
A:=\sup_k\int_0^T\alpha_{\epsilon_k}^t(\mathbb T^2)\,dt<\infty.
\end{equation}
Then $\mu_t(\mathbb T^2)\le M$ for a.e. $t$.
Let $K=\mathbb T^2\times G(2,1)$, and choose a countable dense set
$\{f_j\}_{j\ge1}$ in the unit ball of the real-valued space $C^0(K)$.
For positive measures of mass at most $M$, the metric
\begin{equation}\label{varifold-metric}
d_*(V,W):=\sum_{j=1}^\infty2^{-j}|V(f_j)-W(f_j)|
\end{equation}
metrizes weak convergence.  In particular $d_*(V,W)\le2M$.
Put $d_k(t):=d_*(V_{\epsilon_k}^t,V^t)$.

Fix $L>0$ and a time $t$ at which the conclusions of Steps 2--3 hold.
We claim that
\begin{equation}\label{truncated-varifold-distance}
\mathbf 1_{(0,T)\backslash  \mathcal B_{\epsilon_k,L}}(t)
d_k(t)\longrightarrow0.
\end{equation}
Otherwise, for some $\delta>0$ there would be a subsequence $k_i$ with
$\alpha_{\epsilon_{k_i}}^t(\mathbb T^2)\le L$ and
$d_{k_i}(t)\ge\delta$.  The static compactness theorem
\cite[Theorem 4.1 and its proof]{RogerSchaetzle} gives a further
subsequence whose diffuse varifolds converge to a rectifiable varifold.
Its weight is $\mu_t$, by Step 2.  A rectifiable varifold is uniquely
determined by its weight, since its plane is the approximate tangent of
that weight.  Hence this limit is precisely $V^t$, contradicting
$d_{k_i}(t)\ge\delta$.  This proves
\eqref{truncated-varifold-distance}.  The subsequences in this
contradiction argument may depend on $t$; no common extraction has been
made or assumed.

Dominated convergence, followed by the bad-time estimate, now gives
\begin{equation}\label{varifold-distance-truncation}
\begin{aligned}
\limsup_{k\to\infty}\int_0^T d_k(t)\,dt
&\le \limsup_{k\to\infty}
\int_{(0,T)\setminus\mathcal B_{\epsilon_k,L}}d_k(t)\,dt
+2M\sup_k|\mathcal B_{\epsilon_k,L}|\le\frac{2MA}{L}.
\end{aligned}
\end{equation}
Letting $L\to\infty$ proves
\begin{equation}\label{varifold-L1-time}
\int_0^T d_*(V_{\epsilon_k}^t,V^t)\,dt\longrightarrow0.
\end{equation}
In particular, for every $f\in C(K)$,
\begin{equation}\label{varifold-test-L1-time}
\int_0^T|V_{\epsilon_k}^t(f)-V^t(f)|\,dt\longrightarrow0.
\end{equation}
This follows first for each $f_j$ from \eqref{varifold-metric}, and
then for general $f$ by uniform approximation and the mass bound.
Approximating $F\in C(K\times[0,T])$ uniformly by finite sums
$\sum_i a_i(t)f_i(x,Q)$ gives
\begin{equation}\label{varifold-spacetime-preliminary}
V_{\epsilon_k}^t\,dt\longrightarrow V^t\,dt
\quad\text{as Radon measures on }K\times[0,T].
\end{equation}
This holds along the sequence already fixed, before any additional
extraction.

Finally, by \eqref{varifold-L1-time} we can choose a subsequence $k_i$, independent of $t$, such that
\begin{equation}\label{varifold-summable-extraction}
\int_0^T d_{k_i}(t)\,dt\le2^{-i}.
\end{equation}
Tonelli's theorem gives $\sum_i d_{k_i}(t)<\infty$ for a.e. $t$, so
$V_{\epsilon_{k_i}}^t\to V^t$ for a.e. $t$.  Relabel this single
subsequence as $\epsilon_k$.  All preceding convergences are preserved.

\medskip
\noindent
\emph{Step 6: convergence of the diffuse mean curvature.}
By \eqref{Heps-L2},
 and Lemma \ref{thmB.3} for
measure--function pairs,
there exist a further subsequence and
$\widehat H\in L^2(\mu;\mathbb R^2)$ such that
\begin{equation}\label{liu-revision-equ1}
(\widetilde\mu_{\epsilon_k},H_{\epsilon_k})
\to
(\mu,\widehat H).
\end{equation}
We identify $\widehat H$ with $H_{\mu_t}$.
To this end, let   $X\in C^1(\mathbb T^2\times[0,T];\mathbb R^2)$, and let $V_{\epsilon_k}^t$ be the diffuse varifold associated with
$u_{\epsilon_k}(\cdot,t)$.  Recall from \eqref{RS-first-variation} and \eqref{RS-first-variation+1} that 
$$
\delta V_{\epsilon_k}^t(X)
=
-\int_{\mathbb T^2}
H_{\epsilon_k}\cdot X\,
d(\mu_{\epsilon_k}^t+\xi_{\epsilon_k}^t)
+
R_{\epsilon_k}^t(X),
$$
where
$$
\left|
\int_0^T R_{\epsilon_k}^t(X(\cdot,t))\,dt
\right|
\le
C\|\nabla X\|_{C^0}
|\xi_{\epsilon_k}|(\mathbb T^2\times(0,T))
\to0.
$$
By \eqref{varifold-spacetime-preliminary}, applied to the continuous
Grassmannian test function $(x,Q,t)\mapsto Q:\nabla_xX(x,t)$,
\begin{equation}\label{integrated-first-variation-convergence}
\begin{aligned}
\lim_{k\to\infty}\int_0^T\delta V_{\epsilon_k}^t(X(\cdot,t))\,dt
&=\lim_{k\to\infty}\int_0^T\int Q:\nabla_xX(x,t)
\,dV_{\epsilon_k}^t(x,Q)\,dt\\
&=\int_0^T\delta V^t(X(\cdot,t))\,dt.
\end{aligned}
\end{equation}
Combining this with \eqref{liu-revision-equ1} and
\eqref{discrepancy-total-mass} yields
\begin{equation}\label{integrated-H-identification}
\int_0^T\delta V^t(X(\cdot,t))\,dt
=-\int_0^T\int_{\mathbb T^2}\widehat H\cdot X\,d\mu_t\,dt.
\end{equation}
By Step 3, the left-hand side is also
$-\int_0^T\int H_{\mu_t}\cdot X\,d\mu_t\,dt$.
By uniqueness of the generalized mean curvature,
$$
\widehat H=H_{\mu_t}
\qquad\text{for $\mu_tdt$-a.e.\ $(x,t)$}.
$$
This proves \eqref{H-pair-convergence}.
Notice concretely that
$$
H_{\epsilon_k}\,d\widetilde\mu_{\epsilon_k}
=
v_{\epsilon_k}\nabla u_{\epsilon_k}\,dx\,dt.
$$
Thus \eqref{H-pair-convergence} is equivalently the vector-measure
convergence
$$
v_{\epsilon_k}\nabla u_{\epsilon_k}\,dx\,dt
\to
H_{\mu_t}\,d\mu_t\,dt.
$$

\medskip
\noindent
\emph{Step 7: compactness and normality of the generalized velocity.}
By \eqref{Veps-L2} and \eqref{MR-uniform-bounds},
$$
\sup_k
\int
|\mathcal V_{\epsilon_k}|^2
\,d\widetilde\mu_{\epsilon_k}
<\infty.
$$
Together with \eqref{tilde-mu-limit}, Lemma \ref{thmB.3}  therefore gives a further subsequence and
$\mathcal V\in L^2(\mu;\mathbb R^2)$ such that
$$
(\widetilde\mu_{\epsilon_k},\mathcal V_{\epsilon_k})
\to(\mu,\mathcal V)
$$
and
$$
\int_0^T\int_{\mathbb T^2}
|\mathcal V|^2\,d\mu_t\,dt
\le
\liminf_{k\to\infty}
\int
|\mathcal V_{\epsilon_k}|^2
d\widetilde\mu_{\epsilon_k}.
$$
Together with \eqref{Veps-L2}, this proves
\eqref{velocity-liminf}.

By \eqref{V-pair-convergence}, \eqref{varifold-spacetime-preliminary},
and \eqref{discrepancy-total-mass}, the normality of
$\mathcal V_{\epsilon_k}$ passes to the limit by the argument of
\cite[Lemma 6.3]{MugnaiRoger}; see also
\cite[Proposition 3.2]{Moser}. Hence
$P(x,t)\mathcal V(x,t)=0$ for $\mu$-a.e.\ $(x,t)$,
which proves \eqref{orthogonality1}.

\medskip
\noindent
\emph{Step 8: the $L^2$-flow identity.}
It remains to verify the distributional condition in   Definition \ref{def-L2-flow} of an
$L^2$-flow.  Let
$\eta\in C_c^1(\mathbb T^2\times(0,T))$.
A direct differentiation of the diffuse energy gives
$$
\int_0^T\int_{\mathbb T^2}
\partial_t\eta\,d\mu_\epsilon^t\,dt
=
-\int_0^T\int_{\mathbb T^2}
\eta\,v_\epsilon\,\partial_tu_\epsilon\,dx\,dt
+
\int_0^T\int_{\mathbb T^2}
\epsilon\partial_tu_\epsilon
\nabla u_\epsilon\cdot\nabla\eta\,dx\,dt.
$$
By the definition of $\mathcal V_\epsilon$,
$$
\int_0^T\int_{\mathbb T^2}
\nabla\eta\cdot\mathcal V_\epsilon\,
d\widetilde\mu_\epsilon
=
-\int_0^T\int_{\mathbb T^2}
\epsilon\partial_tu_\epsilon
\nabla u_\epsilon\cdot\nabla\eta\,dx\,dt.
$$
Hence
\begin{equation}\label{approx-L2-flow}
\begin{aligned}
&\int_0^T\int_{\mathbb T^2}
\partial_t\eta\,d\mu_\epsilon^t\,dt
+
\int_0^T\int_{\mathbb T^2}
\nabla\eta\cdot\mathcal V_\epsilon\,
d\widetilde\mu_\epsilon
=
-\int_0^T\int_{\mathbb T^2}
\eta\,v_\epsilon\,\partial_tu_\epsilon\,dx\,dt.
\end{aligned}
\end{equation}
Young's inequality gives
$$
\left|
\int_0^T\int_{\mathbb T^2}
\eta\,v_\epsilon\,\partial_tu_\epsilon\,dx\,dt
\right|
\le
\frac12\|\eta\|_{C^0}
\int_0^T\int_{\mathbb T^2}
\left(
\frac{|v_\epsilon|^2}{\epsilon}
+
\epsilon|\partial_tu_\epsilon|^2
\right)dx\,dt
\le
C(T)\|\eta\|_{C^0}.
$$
Passing to the limit in \eqref{approx-L2-flow}, using
$\mu_{\epsilon_k}\to\mu$ and
\eqref{V-pair-convergence}, yields
$$
\left|
\int_0^T\int_{\mathbb T^2}
\left(
\partial_t\eta+
\nabla\eta\cdot\mathcal V
\right)d\mu_t\,dt
\right|
\le
C(T)\|\eta\|_{C^0}.
$$
Together with the normality of $\mathcal V$, this is precisely the
$L^2$-flow condition.

\medskip
\noindent
\emph{Step 9: Renormalization.}
Set $\widehat\mu_t:=\sigma^{-1}\mu_t$. By definition, the generalized mean curvature of
$\widehat\mu_t$ is the same as that of $\mu_t$. It remains to show that   $\mathcal V$ is a generalized velocity for
$\{\widehat\mu_t\}_{t\in(0,T)}$. Indeed,
$$
\int_0^T\int |\mathcal V|^2\,d\widehat\mu_tdt
=
\sigma^{-1}
\int_0^T\int |\mathcal V|^2\,d\mu_tdt<\infty,
$$
the tangent planes of $\widehat\mu_t$ and $\mu_t$ agree, and
$$
\begin{aligned}
&\left|
\int_0^T\int
(\partial_t\eta+\nabla\eta\cdot\mathcal V)
\,d\widehat\mu_tdt
\right|
\\
&\qquad=
\sigma^{-1}
\left|
\int_0^T\int
(\partial_t\eta+\nabla\eta\cdot\mathcal V)
\,d\mu_tdt
\right|
\le
\sigma^{-1}C_T\|\eta\|_{C^0}.
\end{aligned}
$$
Hence $\{\widehat\mu_t\}_{t\in(0,T)}$ is an $L^2$-flow with generalized
velocity $\mathcal V$.
\end{proof}

\begin{corollary}[Varifold convergence]\label{varifold-convergence-coro}
Let $\epsilon_k\downarrow0$ be the subsequence fixed in Step 5 of
Proposition \ref{MR-compactness} and retained in the subsequent steps.
For a.e. $t\in(0,T)$, let $V_{\epsilon_k}^t$ be the diffuse varifold
in \eqref{def:diffuse-varifold}, and let $V^t$ be the rectifiable
one-dimensional varifold with weight $\mu_t$ and tangent projection
$P(x,t)$.  Then
\begin{equation}\label{varifold-convergence}
V_{\epsilon_k}^t\to V^t
\qquad\text{as varifolds for a.e. }t\in(0,T),
\end{equation}
along this one subsequence, independent of $t$.  Moreover,
\begin{equation}\label{varifold-convergence+1}
V_{\epsilon_k}^t\,dt\to V^t\,dt
\quad\text{as Radon measures on }
\mathbb T^2\times G(2,1)\times(0,T).
\end{equation}
More precisely, \eqref{varifold-L1-time} and
\eqref{varifold-test-L1-time} hold \footnote{these stronger time-integrated
convergences hold even before the additional extraction in
\eqref{varifold-summable-extraction}.}.
\end{corollary}
\begin{proof}
Step 5 proves \eqref{varifold-L1-time} by truncating the diffuse
curvature at a fixed level $L$, using static compactness and uniqueness
of the rectifiable limit on the good times, and then letting
$L\to\infty$.  It proves \eqref{varifold-convergence+1} directly from
this convergence in time.  The summable extraction
\eqref{varifold-summable-extraction} gives
\eqref{varifold-convergence} along a single subsequence.  Thus neither
conclusion relies on integrating convergence along time-dependent
subsequences.
\end{proof}

\section{Fluid velocity compactness via  stream-function}\label{sec:fluid-stream}

Fix $T>0$ throughout this section, and let $\epsilon_k\downarrow0$ be the
subsequence obtained in Proposition \ref{MR-compactness} and in Corollary \ref{varifold-convergence-coro}.  After passing to a
further subsequence, which is not relabeled, we use the same sequence
$\epsilon_k$ for the compactness of the fluid velocities.

\begin{lemma}[Strong compactness of the fluid velocity]\label{fluid-compactness}
There exists
\begin{equation}\label{liu-revision-equ8}
\mathbf v\in L^\infty([0,T];L^2(\mathbb T^2;\mathbb R^2))
\cap L^2([0,T];H^1(\mathbb T^2;\mathbb R^2))
\end{equation}
such that
\begin{equation}\label{liu-revision-equ2}
\mathbf v_{\epsilon_k}\stackrel{*}{\rightharpoonup}\mathbf v
\quad\text{in }L^\infty([0,T];L^2(\mathbb T^2;\mathbb R^2)),
\qquad
\mathbf v_{\epsilon_k}\rightharpoonup\mathbf v
\quad\text{in }L^2([0,T];H^1(\mathbb T^2;\mathbb R^2)).
\end{equation}
Moreover,  for every fixed $s\in (1,\infty)$,
\begin{equation}\label{u-strong}
\mathbf v_{\epsilon_k}\to\mathbf v
\qquad\text{strongly in }
L^2([0,T];L^s(\mathbb T^2;\mathbb R^2)).
\end{equation}
 \end{lemma}

\begin{proof}
The energy inequality \eqref{diffuse-energy} gives a uniform bound for
$\mathbf v_{\epsilon_k}$ in
$L^\infty([0,T];L^2)\cap L^2([0,T];H^1)$.  In two dimensions,
Ladyzhenskaya's inequality gives the uniform bound
\begin{equation}\label{liu-revision-equ3}
\|\mathbf v_{\epsilon_k}\|_{ L^4(\mathbb T^2\times[0,T])}\leq C<\infty.
\end{equation}
Hence
$\mathbf v_{\epsilon_k}\otimes\mathbf v_{\epsilon_k}$ is uniformly bounded
in $L^2(\mathbb T^2\times[0,T])$, and the advection term $(\mathbf v_\epsilon\cdot\nabla)\mathbf v_\epsilon$  is uniformly bounded in
$
L^2([0,T];H^{-1}(\mathbb T^2 )).
$
The capillary stress satisfies, for all $t$,
$$
\|\epsilon_k\nabla u_{\epsilon_k}\otimes\nabla u_{\epsilon_k}\|_{L^1(\mathbb T^2)}
\le 2\mu_{\epsilon_k}^t(\mathbb T^2),
$$
and thus  the right-hand side of equation \eqref{NSAC}  is uniformly bounded in $L^\infty([0,T];W^{-1,1}(\mathbb T^2 ))$ by the
energy inequality.  Using the boundedness of Leray's projector $\mathbb P$ on Sobolev spaces $H^s$ (for any $s\in \mathbb R$) \footnote{
 The proof for the $\mathbb R^n$ case  can be found in \cite[Lemma 1.15]{MR4475666}. The same argument applies on the torus,
since the Leray projection is the Fourier multiplier given by
\[\widehat{\mathbb P f}(k)
=
\begin{cases}
\displaystyle
\left(I-\frac{k\otimes k}{|k|^2}\right)\widehat f(k),
& k\neq0,\\[6pt]
\widehat f(0),& k=0.
\end{cases}
\]}, we deduce that $\partial_t \mathbf v_{\epsilon_k}$ is uniformly bounded in $L^2([0,T];H^{-3}(\mathbb T^2))$. Since $H^1(\mathbb T^2)\hookrightarrow L^s(\mathbb T^2)\hookrightarrow H^{-3}(\mathbb T^2)$,  the Aubin--Lions compactness Theorem gives \eqref{u-strong}. The weak and weak-star
convergences follow from the energy bounds.
\end{proof}

We next transfer the diffuse density estimate \eqref{diffuse-density-est} to the limiting measure.  This
also gives the $L^2$ trace of the limiting fluid velocity that will be used in
Proposition \ref{normal-transport}.

\begin{lemma}\label{limit-fluid-density-trace}
For a.e. $t\in(0,T)$ define
$$
D(t):=
\sup_{x\in\mathbb T^2,\ 0<r\leq 1}
\frac{\mu_t(B_r(x))}{r}.
$$
Then
\begin{equation}\label{limit-density-L2}
\int_0^T D(t)^2\,dt
\le
\liminf_{k\to\infty}
\int_0^T D_{\epsilon_k}(t)^2\,dt
<\infty.
\end{equation}
Consequently,
\begin{equation}\label{limit-fluid-trace}
\int_0^T\int_{\mathbb T^2}
|\mathbf v|^2\,d\mu_t\,dt<\infty.
\end{equation}
\end{lemma}

\begin{proof}
By \eqref{con-mu_t} in Proposition \ref{MR-compactness},  for every $t$, 
$$
\mu_{\epsilon_k}^t\to\mu_t
\qquad\text{as Radon measures on }\mathbb T^2.
$$
For every open ball $B_r(x)$ with $0<r\leq 1$, the
Portmanteau theorem gives
$$
\mu_t(B_r(x))
\le
\liminf_{k\to\infty}\mu_{\epsilon_k}^t(B_r(x))
\le
r\liminf_{k\to\infty}D_{\epsilon_k}(t).
$$
Taking the supremum over $x\in\mathbb T^2$ and $0<r\leq 1$ yields
$$
D(t)\le\liminf_{k\to\infty}D_{\epsilon_k}(t).
$$
Fatou's lemma and \eqref{density-L2} therefore imply
\eqref{limit-density-L2}.

The slice convergence and the diffuse energy bound also give
$$
\operatorname*{ess\,sup}_{0<t<T}\mu_t(\mathbb T^2)<\infty.
$$
Applying Lemma \ref{trace-lemma} to the measure $\mu_t$ and the function
$|\mathbf v(t)|$, we obtain for a.e. $t$
$$
\begin{aligned}
\int_{\mathbb T^2}|\mathbf v|^2\,d\mu_t
&\le
C\,D(t)
\left(
\|\mathbf v(t)\|_{L^2}\|\nabla\mathbf v(t)\|_{L^2}
+
\|\mathbf v(t)\|_{L^2}^2
\right).
\end{aligned}
$$
The right-hand side is integrable in time because
$$
D\in L^2(0,T),
\qquad
\mathbf v\in L^\infty([0,T];L^2(\mathbb T^2;\mathbb R^2))\cap L^2([0,T];H^1(\mathbb T^2;\mathbb R^2)),
$$
and $\mu_t(\mathbb T^2)$ is 
bounded.  This proves
\eqref{limit-fluid-trace}.
\end{proof}

Because $\mathbf v_{\epsilon_k}$ and $\mathbf v$ are divergence free and
have zero spatial mean, there exist unique mean-zero stream functions
$\Psi_{\epsilon_k}$ and $\Psi$ such that
$$
\mathbf v_{\epsilon_k}=\nabla^\perp\Psi_{\epsilon_k},
\qquad
\mathbf v=\nabla^\perp\Psi.
$$
By Poincar\'e's inequality and \eqref{u-strong},
$$
\Psi_{\epsilon_k}\to\Psi
\qquad\text{strongly in }L^2([0,T];H^1(\mathbb T^2)).
$$
Moreover, the energy bounds give
$$
\Psi_{\epsilon_k}
\quad\text{bounded uniformly in}\quad
L^\infty([0,T];H^1(\mathbb T^2))
\cap
L^2([0,T];H^2(\mathbb T^2)).
$$

\begin{lemma}
\label{stream-trace}
We have
\begin{equation}\label{stream-trace-conv}
\int_0^T\int_{\mathbb T^2}
|\Psi_{\epsilon_k}-\Psi|^2\,d\mu_{\epsilon_k}^t\,dt
\to0.
\end{equation}
The same conclusion holds with
$d(\mu_{\epsilon_k}^t+\xi_{\epsilon_k}^t)$ in place of
$d\mu_{\epsilon_k}^t$.
\end{lemma}

\begin{proof}
Set $F_k=\Psi_{\epsilon_k}-\Psi$.  Since both stream functions have zero
mean,
$$
\|F_k(t)\|_{L^2}
+
\|\nabla F_k(t)\|_{L^2}
\le
C\|\mathbf v_{\epsilon_k}(t)-\mathbf v(t)\|_{L^2}.
$$
Lemma \ref{trace-lemma}, applied to $F_k$, therefore gives
$$
\begin{aligned}
\int_{\mathbb T^2}|F_k|^2\,d\mu_{\epsilon_k}^t
&\le
C\,D_{\epsilon_k}(t)
\|\mathbf v_{\epsilon_k}(t)-\mathbf v(t)\|_{L^2}^2.
\end{aligned}
$$
By \eqref{u-strong} and the uniform $L^\infty([0,T];L^2)$ bound of $\mathbf v_{\epsilon_k}$,
interpolation gives
$$
\mathbf v_{\epsilon_k}-\mathbf v
\to0
\qquad\text{strongly in }L^4([0,T];L^2(\mathbb T^2)).
$$
Hence, using the uniform $L^2$ bound for $D_{\epsilon_k}$,
$$
\begin{aligned}
&\int_0^T
D_{\epsilon_k}(t)
\|\mathbf v_{\epsilon_k}(t)-\mathbf v(t)\|_{L^2}^2\,dt
\to0,
\end{aligned}
$$
which proves
\eqref{stream-trace-conv}.
Finally,
$$
0\le \mu_{\epsilon_k}^t+\xi_{\epsilon_k}^t
=
\epsilon_k|\nabla u_{\epsilon_k}|^2\,dx
\le
2\mu_{\epsilon_k}^t,
$$
so the same convergence holds with
$\mu_{\epsilon_k}^t+\xi_{\epsilon_k}^t$ in place of
$\mu_{\epsilon_k}^t$.
\end{proof}

The point of introducing the stream function is revealed in the proof of Proposition \ref{normal-transport} below.  The natural
compactness gives strong convergence of $\mathbf v_{\epsilon_k}$ only in
$L^2_{t,x}$ and does not give strong convergence of
$\nabla\mathbf v_{\epsilon_k}$.  Thus the trace inequality \eqref{trace-est} cannot be applied
directly to $\mathbf v_{\epsilon_k}-\mathbf v$ on the   measures $\mu_{\epsilon_k}^t$.  In
contrast,
$$
\nabla(\Psi_{\epsilon_k}-\Psi)
=
-(\mathbf v_{\epsilon_k}-\mathbf v)^\perp,
$$
so the same strong $L^2$ convergence becomes strong $H^1$ convergence of the
stream functions, which is exactly the regularity required by the trace
estimate.  Proposition \ref{normal-transport} uses the first variation to
transfer the derivative from $\Psi_{\epsilon_k}$ to the varifold.

\begin{proposition}[Identification of the normal transport]
\label{normal-transport}
Let $\epsilon_k\downarrow0$ be the subsequence obtained in
Proposition \ref{MR-compactness} and Corollary \ref{varifold-convergence-coro}. Let
$
P_{\epsilon_k}
=
I-\nu_{\epsilon_k}\otimes\nu_{\epsilon_k}
$
be the tangent projection in \eqref{projection-eps}, and let $P=P(x,t)$
denote the tangent projection of the limiting integral varifold
$\sigma^{-1}\mu_t$. We set
\begin{equation}
P_{\epsilon_k}^{\perp}:=I-P_{\epsilon_k},
\qquad
P^\perp:=I-P.
\end{equation}
Then, for every
$\zeta\in C^1(\mathbb T^2\times[0,T];\mathbb R^2)$,
\begin{equation}\label{normal-transport-conv}
\begin{aligned}
&\lim_{k\to\infty}
\int_0^T\int_{\mathbb T^2}
P_{\epsilon_k}^\perp\mathbf v_{\epsilon_k}\cdot\zeta\,
d(\mu_{\epsilon_k}^t+\xi_{\epsilon_k}^t)\,dt
=
\int_0^T\int_{\mathbb T^2}
P^\perp\mathbf v\cdot\zeta\,d\mu_t\,dt.
\end{aligned}
\end{equation}
\end{proposition}

\begin{proof}

For a vector $a=(a_1,a_2)$, we write
$
a^\perp=(-a_2,a_1).
$
Notice that this is the vector by  rotating $a$ counterclockwise through an angle of   $\pi/2$, which is consistent with \eqref{liu-revision-equ7}, whereas
$P^\perp=I-P$ denotes the orthogonal projection onto the normal line.

\medskip
\noindent
\emph{Step 1: an identity involving the first variation.}
For every rank-one orthogonal projection $P$ in $\mathbb R^2$, we have  $
J P J^T=P^{\perp}$ where $
J=\left(\begin{array}{cc}
0 & -1 \\
1 & 0
\end{array}\right) 
$.
\footnote{Indeed, 
 if  $P=\tau \otimes \tau$, where $|\tau|=1$,  then $J\tau$ is the unit vector perpendicular to $\tau$. Hence
$
J P J^T=J(\tau \otimes \tau) J^T=(J \tau) \otimes(J \tau) .
$
Since $\{\tau, J \tau\}$ is an orthonormal basis of $\mathbb{R}^2$,
$
I=\tau \otimes \tau+(J \tau) \otimes(J \tau) .
$
 } Let $\tau$ be a unit tangent vector and set $\nu=\tau^\perp$.
Then $P=\tau\otimes\tau$ and $P^\perp=\nu\otimes\nu$, so
\begin{equation}
\label{perp-formula1}
P^\perp\nabla^\perp\Psi
=(\nu\cdot\nabla^\perp\Psi)\nu
=(\tau\cdot\nabla\Psi)\tau^\perp
=(P\nabla\Psi)^\perp.
\end{equation}
Therefore
$$
P_{\epsilon_k}^\perp\mathbf v_{\epsilon_k}\cdot\zeta
=
-P_{\epsilon_k}\nabla\Psi_{\epsilon_k}\cdot\zeta^\perp.
$$
Applying the definition 
$$
(\delta V_{\epsilon_k}^t)(X)
=
\int_{\mathbb T^2}
P_{\epsilon_k}:\nabla X\,d\mu_{\epsilon_k}^t,
$$
with
$X=-\Psi_{\epsilon_k}\zeta^\perp$, we obtain
\begin{equation}\label{stream-first-variation}
\begin{aligned}
&\int_{\mathbb T^2}
P_{\epsilon_k}^\perp\mathbf v_{\epsilon_k}\cdot\zeta\,
d\mu_{\epsilon_k}^t
 =
(\delta V_{\epsilon_k}^t)
(-\Psi_{\epsilon_k}\zeta^\perp)
+
\int_{\mathbb T^2}
\Psi_{\epsilon_k}
P_{\epsilon_k}:\nabla \zeta^\perp\,
d\mu_{\epsilon_k}^t.
\end{aligned}
\end{equation}
By the   first-variation identity
\eqref{RS-first-variation+1},
$$
\begin{aligned}
(\delta V_{\epsilon_k}^t)(X)
={}&
-\int_{\mathbb T^2}
H_{\epsilon_k}\cdot X\,
d(\mu_{\epsilon_k}^t+\xi_{\epsilon_k}^t)+
\int_{\mathbb T^2}
(\nu_{\epsilon_k}\otimes\nu_{\epsilon_k}):\nabla X\,
d\xi_{\epsilon_k}^t.
\end{aligned}
$$
Hence, after integration in time,
\begin{equation}\label{stream-identity-2}
\begin{aligned}
\int_0^T\int_{\mathbb T^2}
P_{\epsilon_k}^\perp\mathbf v_{\epsilon_k}\cdot\zeta\,
d\mu_{\epsilon_k}^t\,dt
=&
\int_0^T\int_{\mathbb T^2}
\Psi_{\epsilon_k}
H_{\epsilon_k}\cdot\zeta^\perp\,
d(\mu_{\epsilon_k}^t+\xi_{\epsilon_k}^t)\,dt
\\
&\quad+
\int_0^T\int_{\mathbb T^2}
\Psi_{\epsilon_k}
P_{\epsilon_k}:\nabla \zeta^\perp\,
d\mu_{\epsilon_k}^t\,dt
+
R_k,
\end{aligned}
\end{equation}
where
$$
R_k
=
-\int_0^T\int_{\mathbb T^2}
(\nu_{\epsilon_k}\otimes\nu_{\epsilon_k})
:
\nabla(\Psi_{\epsilon_k}\zeta^\perp)
\,d\xi_{\epsilon_k}^t\,dt.
$$

\medskip
\noindent
\emph{Step 2: the discrepancy term vanishes.}
Since
$$
\nabla(\Psi_{\epsilon_k}\zeta^\perp)
=
\zeta^\perp\otimes\nabla\Psi_{\epsilon_k}
+
\Psi_{\epsilon_k}\nabla\zeta^\perp
$$
up to the harmless convention for the order of tensor indices, and
$
|\nabla\Psi_{\epsilon_k}|=|\mathbf v_{\epsilon_k}|,
$
we have
\begin{equation}\label{Rk-estimate}
|R_k|
\le
C_\zeta
\int_0^T\int_{\mathbb T^2}
\left(
|\mathbf v_{\epsilon_k}|+|\Psi_{\epsilon_k}|
\right)
\,d|\xi_{\epsilon_k}^t|\,dt.
\end{equation}
Pointwise,
$
|\xi_{\epsilon_k}^t|\le\mu_{\epsilon_k}^t.
$
Therefore
$$
\begin{aligned}
\int |\mathbf v_{\epsilon_k}|\,d|\xi_{\epsilon_k}|
&\le
\left(
\int |\mathbf v_{\epsilon_k}|^2\,d\mu_{\epsilon_k}
\right)^{1/2}
|\xi_{\epsilon_k}|
(\mathbb T^2\times(0,T))^{1/2},
\\
\int |\Psi_{\epsilon_k}|\,d|\xi_{\epsilon_k}|
&\le
\left(
\int |\Psi_{\epsilon_k}|^2\,d\mu_{\epsilon_k}
\right)^{1/2}
|\xi_{\epsilon_k}|
(\mathbb T^2\times(0,T))^{1/2}.
\end{aligned}
$$
The first $L^2$ factor is uniformly bounded by Proposition
\ref{diffuse-trace}. For the second one, the mean-zero condition and
Poincar\'e's inequality give
$$
\|\Psi_{\epsilon_k}(t)\|_{H^1}
\le C\|\mathbf v_{\epsilon_k}(t)\|_{L^2},
$$
and the same diffuse trace estimate gives
$$
\sup_k
\int_0^T\int_{\mathbb T^2}
|\Psi_{\epsilon_k}|^2\,d\mu_{\epsilon_k}^t\,dt<\infty.
$$
Since
$
|\xi_{\epsilon_k}|
(\mathbb T^2\times(0,T))\to0
$
by \eqref{discrepancy-vanish}, \eqref{Rk-estimate} yields
\begin{equation}\label{Rk-zero}
R_k\to0.
\end{equation}

\medskip
\noindent
\emph{Step 3: passage to the limit in the two principal terms.}
We first claim
\begin{equation}\label{H-stream-limit}
\begin{aligned}
&\int_0^T\int_{\mathbb T^2}
\Psi_{\epsilon_k}H_{\epsilon_k}\cdot\zeta^\perp\,
d(\mu_{\epsilon_k}^t+\xi_{\epsilon_k}^t)\,dt
 \longrightarrow
\int_0^T\int_{\mathbb T^2}
\Psi H_{\mu_t}\cdot\zeta^\perp\,d\mu_t\,dt.
\end{aligned}
\end{equation}
Indeed, by Cauchy--Schwarz,
$$
\begin{aligned}
&\left|
\int
(\Psi_{\epsilon_k}-\Psi)
H_{\epsilon_k}\cdot\zeta^\perp\,
d(\mu_{\epsilon_k}+\xi_{\epsilon_k})
\right|
\\
&\qquad\le
\|\zeta\|_{C^0}
\left(
\int
|\Psi_{\epsilon_k}-\Psi|^2
\,d(\mu_{\epsilon_k}+\xi_{\epsilon_k})
\right)^{1/2}
\left(
\int
|H_{\epsilon_k}|^2
\,d(\mu_{\epsilon_k}+\xi_{\epsilon_k})
\right)^{1/2},
\end{aligned}
$$
which tends to zero by Lemma \ref{stream-trace} and \eqref{curvature-liminf}.
To treat  the remaining term containing $\Psi$, we first deduce from  \eqref{liu-revision-equ8} and interpolation that   $\Psi\in L^4([0,T];H^1(\mathbb T^2))$. So  choose
$\Psi^{(m)}\in C^\infty(\mathbb T^2\times[0,T])$ such that
\begin{equation}
\Psi^{(m)}\to\Psi
\qquad\text{in }L^4([0,T];H^1(\mathbb T^2)).
\end{equation}
The diffuse and limiting   trace estimate imply 
\begin{align}
\lim_{m\to\infty}
\sup_k
\int_0^T\int_{\mathbb T^2}
|\Psi^{(m)}-\Psi|^2\,d\mu_{\epsilon_k}^t\,dt
=0,\\
\lim_{m\to\infty}
\int_0^T\int_{\mathbb T^2}
|\Psi^{(m)}-\Psi|^2\,d\mu_t\,dt
=0.
\end{align}
For fixed $m$, the measure--function convergence
\eqref{H-pair-convergence} applies to the smooth test field
$\Psi^{(m)}\zeta^\perp$.  Letting first $k\to\infty$ and then
$m\to\infty$ proves \eqref{H-stream-limit}.
Likewise,
\begin{equation}\label{P-stream-limit}
\begin{aligned}
&\int_0^T\int_{\mathbb T^2}
\Psi_{\epsilon_k}
P_{\epsilon_k}:\nabla\zeta^\perp\,
d\mu_{\epsilon_k}^t\,dt
\longrightarrow
\int_0^T\int_{\mathbb T^2}
\Psi P:\nabla\zeta^\perp\,d\mu_t\,dt.
\end{aligned}
\end{equation}
To see this, first replace $\Psi_{\epsilon_k}$ by $\Psi$ using
Lemma \ref{stream-trace}. Then approximate $\Psi$ by the smooth functions
$\Psi^{(m)}$ above and use varifold convergence \eqref{varifold-convergence+1} for fixed $m$.

Combining \eqref{stream-identity-2}, \eqref{Rk-zero},
\eqref{H-stream-limit}, and \eqref{P-stream-limit}, we obtain
\begin{equation}\label{stream-limit-expression}
\begin{aligned}
&\lim_{k\to\infty}
\int_0^T\int_{\mathbb T^2}
P_{\epsilon_k}^\perp\mathbf v_{\epsilon_k}\cdot\zeta\,
d\mu_{\epsilon_k}^t\,dt
\\
&=
\int_0^T\int_{\mathbb T^2}
\Psi H_{\mu_t}\cdot\zeta^\perp\,d\mu_t\,dt
+
\int_0^T\int_{\mathbb T^2}
\Psi P:\nabla\zeta^\perp\,d\mu_t\,dt.
\end{aligned}
\end{equation}

\medskip
\noindent
\emph{Step 4: identification of the limiting expression.}
By \eqref{limit-fluid-trace},
$
\mathbf v\in L^2(d\mu_tdt;\mathbb R^2).
$
We claim that
\begin{equation}\label{limit-stream-ibp}
\begin{aligned}
\text{R.H.S of } \eqref{stream-limit-expression} =\int_0^T\int_{\mathbb T^2}
P^\perp\mathbf v\cdot\zeta\,d\mu_t\,dt.
\end{aligned}
\end{equation}

To justify this identity, approximate the fixed limiting velocity
$\mathbf v$ by smooth divergence-free vector fields
$\mathbf v^{(m)}$ such that $\|\mathbf v^{(m)}\|_{L^\infty_tL^2_x}\leq C$ and 
$$
\mathbf v^{(m)}\to\mathbf v
\qquad\text{strongly in }
L^2([0,T];H^1(\mathbb T^2))
$$
and, after passage to a subsequence,
$$
\|\mathbf v^{(m)}(t)-\mathbf v(t)\|_{L^2}\to0
\qquad\text{for a.e.\ }t.
$$
Let $\Psi^{(m)}$ be their mean-zero stream functions.

Thanks to \eqref{limit-density-L2}, the Meyers--Ziemer estimate gives
\begin{equation}
\label{mollif-v}
\int_0^T\int_{\mathbb T^2}
\left(|\Psi^{(m)}-\Psi|^2+|\mathbf v^{(m)}-\mathbf v|^2\right)\,d\mu_t\,dt\to0.
\end{equation}
Indeed, with
$$
a_m(t)=
\|\mathbf v^{(m)}(t)-\mathbf v(t)\|_{L^2},
\qquad
b_m(t)=
\|\nabla(\mathbf v^{(m)}-\mathbf v)(t)\|_{L^2},
$$
the trace estimate gives a bound by
$$
C\int_0^T
D(t)
\big(a_m(t)b_m(t)+a_m(t)^2\big)\,dt.
$$
Using
$
D(t)\,a_mb_m
\le
\frac12D^2(t)a_m^2+\frac12b_m^2,
$
the first term tends to zero by dominated convergence, while
$b_m\to0$ in $L^2(0,T)$.

For every smooth $\Psi^{(m)}$, the first-variation and \eqref{perp-formula1} gives
$$
\begin{aligned}
\int
P^\perp\mathbf v^{(m)}\cdot\zeta\,d\mu
&=
(\delta V)(-\Psi^{(m)}\zeta^\perp)
+
\int
\Psi^{(m)}P:\nabla\zeta^\perp\,d\mu
\\
&=
\int
\Psi^{(m)}H_{\mu_t}\cdot\zeta^\perp\,d\mu
+
\int
\Psi^{(m)}P:\nabla\zeta^\perp\,d\mu\\
&=\int 
P^\perp\mathbf v^{(m)}\cdot\zeta\,d\mu.
\end{aligned}
$$
Letting $m\to\infty$ and using \eqref{mollif-v} and the method in \emph{Step 3}, we get  \eqref{limit-stream-ibp}.  Together with
\eqref{stream-limit-expression}, this gives
\begin{equation}\label{normal-mu-limit}
\begin{aligned}
&\lim_{k\to\infty}
\int_0^T\int_{\mathbb T^2}
P_{\epsilon_k}^\perp\mathbf v_{\epsilon_k}\cdot\zeta\,
d\mu_{\epsilon_k}^t\,dt
=
\int_0^T\int_{\mathbb T^2}
P^\perp\mathbf v\cdot\zeta\,d\mu_t\,dt.
\end{aligned}
\end{equation}

\medskip
\noindent
\emph{Step 5: replacing $\mu_{\epsilon_k}$ by
$\mu_{\epsilon_k}+\xi_{\epsilon_k}$.}
Finally,
$$
\begin{aligned}
&\left|
\int_0^T\int_{\mathbb T^2}
P_{\epsilon_k}^\perp\mathbf v_{\epsilon_k}\cdot\zeta\,
d\xi_{\epsilon_k}^t\,dt
\right|
\\
&\qquad\le
\|\zeta\|_{C^0}
\left(
\int_0^T\int_{\mathbb T^2}
|\mathbf v_{\epsilon_k}|^2\,d\mu_{\epsilon_k}^t\,dt
\right)^{1/2}
|\xi_{\epsilon_k}|
(\mathbb T^2\times(0,T))^{1/2}.
\end{aligned}
$$
The first factor is uniformly bounded by Proposition
\ref{diffuse-trace}, whereas the second tends to zero by
\eqref{discrepancy-vanish}. Hence
$$
\int
P_{\epsilon_k}^\perp\mathbf v_{\epsilon_k}\cdot\zeta\,
d\xi_{\epsilon_k}\to0.
$$
Combining this with \eqref{normal-mu-limit} proves
\eqref{normal-transport-conv}.
\end{proof}

\begin{corollary}[Limiting motion law]\label{motion-law}
The generalized velocity in Proposition \ref{MR-compactness} satisfies
\begin{equation}\label{velocity-law}
\mathcal V=H_{\mu_t}+P^\perp\mathbf v
\qquad\text{for $\mu_t\,dt$-a.e.\ $(x,t)$}.
\end{equation}
\end{corollary}

\begin{proof}
Pass to the limit in \eqref{exact-diffuse-velocity}, using Proposition
\ref{MR-compactness} and Proposition \ref{normal-transport}.
\end{proof}


\section{Proof of Brakke's inequality  \eqref{brakke-inequ}}\label{section-brakke}

\begin{proposition}[Localized coupled energy inequality]\label{loc-ene-ineq}
Let $\epsilon_k\downarrow0$ be the subsequence obtained in
Proposition \ref{MR-compactness} and $\mathbf v$ be the limit in Proposition \ref{fluid-compactness}.
Set $\lambda_t=\sigma^{-1}\mu_t$. There exists a pressure  
\begin{equation}\label{liu-revision-equ4}
\Pi \in L_{\text {loc }}^2 ([0, \infty) ; L^r (\mathbb{T}^2 ) ) \quad \text { for every } r\in (1,2)
\end{equation}
such that    for every nonnegative
$\phi\in C_c^\infty(\mathbb T^2\times [0,\infty))$,
\begin{equation}\label{localized-coupled-energy}
\begin{aligned}
&\int_0^\infty \int_{\mathbb T^2}
\phi|\nabla\mathbf v|^2\,dx\,dt
+\int_0^\infty \int_{\mathbb T^2}
\phi|H_{\mu_t}|^2\,d\lambda_t\,dt
\\
&\quad\le
\int_0^\infty \int_{\mathbb T^2}
\frac{1}{2}|\mathbf v|^2
(\partial_t\phi+\Delta\phi)\,dx\,dt
\\
&\qquad+
\int_0^\infty \int_{\mathbb T^2}
\left(\frac{1}{2}|\mathbf v|^2+\Pi\right)
\mathbf v\cdot\nabla\phi\,dx\,dt
\\
&\qquad+
\int_0^\infty \int_{\mathbb T^2}
\left[
\partial_t\phi+
(H_{\mu_t}+P^\perp\mathbf v)\cdot\nabla\phi
\right]\,d\lambda_t\,dt\\
&\qquad+
\frac12\int_{\mathbb T^2}
\phi(x,0)|\mathbf v_0(x)|^2\,dx
+\int_{\mathbb T^2}
\phi(x,0)\,d\|\nabla\mathcal X_{E_0}\|(x).
\end{aligned}
\end{equation}
\end{proposition}

\begin{proof}
We shall abbreviate $\int_0^\infty \int_{\mathbb T^2}$ by $\iint$.
Write
\begin{equation}\label{localized-diffuse-notation}
k_\epsilon=\frac12|\mathbf v_\epsilon|^2,
\qquad
e_\epsilon=\frac{\epsilon}{2}|\nabla u_\epsilon|^2
+\frac1\epsilon W(u_\epsilon),
\qquad
\pi_\epsilon=\Pi_\epsilon+\sigma^{-1}e_\epsilon.
\end{equation}
The pressures are defined up to functions of time, and we normalize
$\pi_\epsilon$ to have zero spatial mean. Since
\begin{equation}\label{localized-capillary-identity}
v_\epsilon\nabla u_\epsilon
=\nabla e_\epsilon-
\operatorname{div}
(\epsilon\nabla u_\epsilon\otimes\nabla u_\epsilon),
\end{equation}
the momentum equation in \eqref{NSAC} becomes
\begin{equation}\label{localized-modified-NS}
\partial_t\mathbf v_\epsilon
+\operatorname{div}
(\mathbf v_\epsilon\otimes\mathbf v_\epsilon)
-\Delta\mathbf v_\epsilon+\nabla\pi_\epsilon
=\sigma^{-1}v_\epsilon\nabla u_\epsilon.
\end{equation}
Taking its scalar product with $\mathbf v_\epsilon$ gives
\begin{equation}\label{localized-kinetic-identity}
\begin{aligned}
\partial_tk_\epsilon+|\nabla\mathbf v_\epsilon|^2
&=\Delta k_\epsilon
-\operatorname{div}
\big((k_\epsilon+\pi_\epsilon)\mathbf v_\epsilon\big)
+
\sigma^{-1}v_\epsilon
\mathbf v_\epsilon\cdot\nabla u_\epsilon.
\end{aligned}
\end{equation}
On the other hand, \eqref{stationary} and \eqref{AC-v} yield
\begin{equation}\label{localized-phase-identity}
\partial_te_\epsilon+\frac{|v_\epsilon|^2}{\epsilon}
=
\operatorname{div}
(\epsilon\partial_tu_\epsilon\nabla u_\epsilon)
-v_\epsilon\mathbf v_\epsilon\cdot\nabla u_\epsilon.
\end{equation}
Adding $\sigma^{-1}$ times \eqref{localized-phase-identity} to
\eqref{localized-kinetic-identity} cancels the mixed term exactly
and gives
\begin{equation}\label{localized-coupled-identity}
\begin{aligned}
&\partial_t\left(k_\epsilon+\sigma^{-1}e_\epsilon\right)
+|\nabla\mathbf v_\epsilon|^2
+\frac{|v_\epsilon|^2}{\sigma\epsilon}
\\
&\qquad=
\Delta k_\epsilon
-\operatorname{div}\bigl((k_\epsilon+\pi_\epsilon)
\mathbf v_\epsilon\bigr)
+\sigma^{-1}\operatorname{div}
\bigl(\epsilon\partial_tu_\epsilon\nabla u_\epsilon\bigr).
\end{aligned}
\end{equation}
Multiplication by $\phi$ and integration by parts therefore give
\begin{equation}\label{localized-diffuse-energy}
\begin{aligned}
&\iint
\phi|\nabla\mathbf v_\epsilon|^2\,dx\,dt
+\frac1\sigma\iint
\phi\frac{|v_\epsilon|^2}{\epsilon}\,dx\,dt
\\
&\quad=
\iint
k_\epsilon(\partial_t\phi+\Delta\phi)\,dx\,dt
+\frac1\sigma\iint
e_\epsilon\partial_t\phi\,dx\,dt
\\
&\qquad+
\iint
(k_\epsilon+\pi_\epsilon)
\mathbf v_\epsilon\cdot\nabla\phi\,dx\,dt
\\
&\qquad-
\frac1\sigma\iint
\epsilon\partial_tu_\epsilon
\nabla u_\epsilon\cdot\nabla\phi\,dx\,dt\\
&\qquad +\int_{\mathbb T^2}
\phi(x,0)\left(
\frac12|\mathbf v_{\epsilon,0}(x)|^2
+\sigma^{-1}e_\epsilon(x,0)
\right)\,dx. 
\end{aligned}
\end{equation}
We now pass to the limit along $\epsilon=\epsilon_k$.
By  \eqref{liu-revision-equ2}, \eqref{u-strong} and  interpolation, we get 
\begin{equation}\label{localized-fluid-strong}
\begin{aligned}
\mathbf v_{\epsilon_k}&\to\mathbf v
&&\text{strongly in }L^2(0,T;L^s(\mathbb T^2)),
\quad 2\le s<\infty,\\
\mathbf v_{\epsilon_k}&\to\mathbf v
&&\text{strongly in }L^p(\mathbb T^2\times(0,T)),
\quad 2\le p<4.
\end{aligned}
\end{equation}
The case $p=3$ gives convergence of the cubic kinetic flux in \eqref{localized-diffuse-energy}, while
the case $p=2$ gives convergence of the term containing
$k_{\epsilon_k}(\partial_t\phi+\Delta\phi)$.

To treat the terms involving    $\pi_\epsilon$, set
$\mathbf f_\epsilon=\sigma^{-1}v_\epsilon\nabla u_\epsilon$.
Cauchy--Schwarz and \eqref{diffuse-energy} imply
\begin{equation}\label{localized-force-bound}
\begin{aligned}
\|\mathbf f_\epsilon(t)\|_{L^1}^2
&\le \frac1{\sigma^2}
\left(\int_{\mathbb T^2}
\frac{|v_\epsilon|^2}{\epsilon}\,dx\right)
\left(\int_{\mathbb T^2}
\epsilon|\nabla u_\epsilon|^2\,dx\right),\qquad
\int \|\mathbf f_\epsilon(t)\|_{L^1}^2\,dt
&\le C.
\end{aligned}
\end{equation}
Taking the divergence of \eqref{localized-modified-NS} and using
$\operatorname{div}\mathbf v_\epsilon=0$, we obtain
\begin{equation}\label{localized-pressure-equation}
\Delta\pi_\epsilon
=-\sum\nolimits_{i,j=1}^2
\partial_i\partial_j
(v_{\epsilon,i}v_{\epsilon,j})
+\sigma^{-1}\operatorname{div}
(v_\epsilon\nabla u_\epsilon).
\end{equation}
Decompose the mean-zero pressure as
$\pi_\epsilon=\pi_\epsilon^{N}+\pi_\epsilon^{C}$, where both components have zero spatial mean and 
\begin{equation}\label{localized-pressure-decomposition}
\Delta\pi_\epsilon^{N}
=-\sum\nolimits_{i,j=1}^2
\partial_i\partial_j(v_{\epsilon,i}v_{\epsilon,j}),
\qquad
\Delta\pi_\epsilon^{C}
=\operatorname{div}\mathbf f_\epsilon.
\end{equation}
 The uniform $L^4$ bound for $\mathbf v_\epsilon$ in  \eqref{liu-revision-equ3} and the
Calderón--Zygmund estimate give
\begin{equation}\label{localized-NS-pressure-bound}
\|\pi_\epsilon^{N}\|_{L^2(\mathbb T^2\times(0,T))}
\le C.
\end{equation}
\footnote{Indeed, for each nonzero Fourier mode $k\in\mathbb (2Z)^2$,
$
\widehat{\pi_\epsilon^{N}}(k)
=-\sum_{i,j=1}^2\frac{k_i k_j}{|k|^2}
\widehat{v_{\epsilon,i}v_{\epsilon,j}}(k),
$
so the estimate follows directly from Parseval's identity.
}
Moreover, the periodic kernel $\mathcal K=\nabla G$, where $G$
is the mean-zero Green function of $\Delta$ on $\mathbb T^2$,
satisfies
\[
|\mathcal K(x)|
\le C\left(1+d_{\mathbb T^2}(x,0)^{-1}\right).
\]
Hence $\mathcal K\in L^r(\mathbb T^2)$ for every $r\in (1,2)$.
So  Young's convolution inequality  
gives
\[
\|\pi_\epsilon^C(t)\|_{L^r}
\leq \|\mathcal{K}*\mathbf f_\epsilon(t)\|_{L^r}
\le \|\mathcal{K}\|_{L^r}\|\mathbf f_\epsilon(t)\|_{L^1}
\]
for $r\in (1,2)$.
  Combined with  \eqref{localized-force-bound}, this   yields
\begin{equation}\label{localized-capillary-pressure-bound}
\|\pi_\epsilon^{C}\|_{L^2(0,T;L^r(\mathbb T^2))}
\le C_r,
\qquad r\in (1,2).
\end{equation}
Fix $r\in (1,2)$.
After extracting a subsequence, 
\begin{equation}\label{liu-revision-equ6}
\pi^{N}_{\epsilon_k} {\rightharpoonup}\pi^{N}
\quad\text{in }L^2([0,T];L^2(\mathbb T^2)),
\qquad
\pi^C_{\epsilon_k}\rightharpoonup \pi^C
\quad\text{in }L^2([0,T];L^r(\mathbb T^2))
\end{equation}
with $r\in (1,2)$.
Then the pressure
\begin{equation}
    \Pi=\pi^{N}+\pi^{C}
\end{equation} satisfies   \eqref{liu-revision-equ4}.
 Moreover, using 
\eqref{localized-fluid-strong} with $s=r/(r-1)$ for
$\pi_{\epsilon_k}^{C}$, we obtain
\begin{equation}\label{localized-pressure-flux-limit}
\iint
\pi_{\epsilon_k}\mathbf v_{\epsilon_k}\cdot\nabla\phi\,dx\,dt
\longrightarrow
\iint
\pi\mathbf v\cdot\nabla\phi\,dx\,dt.
\end{equation}

The space--time measure convergence accompanying
\eqref{con-mu_t} gives
\begin{equation}\label{localized-phase-energy-limit}
\iint
e_{\epsilon_k}\partial_t\phi\,dx\,dt
\longrightarrow
\iint
\partial_t\phi\,d\mu_t\,dt.
\end{equation}
By \eqref{def:diffusemcf} and \eqref{def: mu-tilde},
$-\epsilon\partial_tu_\epsilon\nabla u_\epsilon\,dx\,dt
=\mathcal V_\epsilon\,d\widetilde\mu_\epsilon.$
Thus \eqref{V-pair-convergence} yields
\begin{equation}\label{localized-phase-flux-limit}
-\iint
\epsilon_k\partial_tu_{\epsilon_k}
\nabla u_{\epsilon_k}\cdot\nabla\phi\,dx\,dt
\longrightarrow
\iint
\mathcal V\cdot\nabla\phi\,d\mu_t\,dt.
\end{equation}
In view of \eqref{velocity-law}, this gives the second to last   term   of \eqref{localized-coupled-energy}.

Finally, weak lower semicontinuity gives
\begin{equation}\label{localized-viscous-liminf}
\iint
\phi|\nabla\mathbf v|^2\,dx\,dt
\le
\liminf_{k\to\infty}
\iint
\phi|\nabla\mathbf v_{\epsilon_k}|^2\,dx\,dt.
\end{equation}
The weighted lower semicontinuity of measure--function pairs,
applied to \eqref{H-pair-convergence}, gives the localized version
of \eqref{curvature-liminf}:
\begin{equation}\label{localized-curvature-liminf}
\begin{aligned}
\iint
\phi|H_{\mu_t}|^2\,d\mu_t\,dt
&\le
\liminf_{k\to\infty}
\iint
\phi\frac{|v_{\epsilon_k}|^2}{\epsilon_k}\,dx\,dt.
\end{aligned}
\end{equation}
Passing to the limit in \eqref{localized-diffuse-energy}, and using
the motion identity \eqref{velocity-law}, proves
\eqref{localized-coupled-energy}.
\end{proof}

\begin{proof}[Proof of \eqref{weak-NS}]
 By \eqref{gradient-measure} and \eqref{discrepancy-vanish},
$$
\epsilon_k|\nabla u_{\epsilon_k}|^2\,dx\,dt
=(\mu_{\epsilon_k}+\xi_{\epsilon_k})dt\to\mu
$$
 as Radon measures on 
$\mathbb T^2\times[0,T].$
Moreover, in \eqref{def:diffuse-varifold}, taking  $$
f(x, Q)=A(x, t):(I-Q)
$$
for continuous matrix-valued $A$ and $Q\in G(2,1)$, and using  \eqref{varifold-convergence+1}, we find
\[
\int A: P_{\epsilon_k}^{\perp} d \mu_{\epsilon_k} \longrightarrow \int A: P^{\perp} d V
\]
where $dV=dV_t dt$ and $V_t$ is the integral varifold with weight $\mu_t$.
This implies that 
$$
\epsilon_k\nabla u_{\epsilon_k}\otimes\nabla u_{\epsilon_k}\,dx\,dt
\to P^\perp\,dV.
$$
Let $\eta\in C_c^1(\mathbb T^2\times[0,\infty);\mathbb R^2)$ .  Then by 
$
P^\perp+P=I_2,$
we have 
$$
 \int P^\perp:\nabla \eta\,dV 
= \int {\rm div }\, \eta\, d\mu_t-  \delta V(\eta)
= \int {\rm div }\, \eta\, d\mu_t+\int H_{\mu_t}\cdot\eta\,d\mu_t.
$$
This gives the convergence of  the capillary term of the $\epsilon$-Navier--Stokes equation (cf. \eqref{NSAC})
\begin{equation}
\iint    \epsilon\nabla u_{\epsilon_k}\otimes\nabla u_{\epsilon_k}:\nabla\eta\,dx\,dt\longrightarrow
 \iint   [H_{\mu_t}\cdot\eta+{\rm div }\, \eta]\,d\mu_t\,dt.
\end{equation}
The convection term passes to the limit by
\eqref{u-strong}, while the viscous term passes by weak convergence in
$L^2_tH^1_x$. Concerning the pressure converges, recall from \eqref{localized-diffuse-notation} that $\Pi_\epsilon=\pi_\epsilon-\sigma^{-1}e_\epsilon.$ We deduce from  \eqref{liu-revision-equ6} and \eqref{con-mu_t} that
\[\int \Pi_{\epsilon_k}\, {\rm div}\,\eta\, dxdt\longrightarrow \int \Pi\, {\rm div}\,\eta\, dxdt-\sigma^{-1}\int {\rm div}\,\eta\,d\mu_t\,dt.
\]
Altogether, we are able to pass to the sequential limit in the weak formulation of \eqref{NSAC} and obtain:
\begin{equation}\label{weak-NS+without-p}
\begin{aligned}
&-\int_0^\infty\int_{\mathbb T^2}
\Big(\mathbf v\cdot\partial_t\eta
+(\mathbf v\otimes\mathbf v):\nabla\eta-\nabla\mathbf v:\nabla\eta\Big)\,
dx\,dt\\&
+\int_0^\infty\int_{\mathbb T^2}{\rm div }\, \eta\, d(\sigma^{-1}\mu_t)\, dt-\int_0^\infty\int_{\mathbb T^2}
 \Pi\,{\rm div}\,\eta\,dx\,dt\\&
\quad =
\int_{\mathbb T^2}\mathbf v_0\cdot\eta(\cdot,0)\,dx
+\int_0^\infty\int_{\mathbb T^2}(
H_{\mu_t}\cdot\eta+ {\rm div }\, \eta)\, d(\sigma^{-1}\mu_t)\, dt.
\end{aligned}
\end{equation}
This is exactly  \eqref{weak-NS}. 
\end{proof}

We next prove the following equality utilizing \eqref{weak-NS}.
\begin{proposition}
For the limiting quantities and for a.e.~$0<t_1<t_2<\infty$, we have
    \begin{equation}\label{weNS}
\begin{split}
  &\left.\frac12\int_{\mathbb T^2} \zeta|{\bf v}|^2\,dx\right|_{t=t_1}^{t_2}  
  +\int_{t_1}^{t_2}\zeta|\nabla{\bf v}|^2\,dxdt
  \\&=\int_{t_1}^{t_2}\int_{\mathbb T^2}\frac12|{\bf v}|^2(\partial_t\zeta+\Delta\zeta)+\left(\frac12|{\bf v}|^2+\Pi\right){\bf v}\cdot\nabla\zeta\,dxdt 
\\&\quad +\int_{t_1}^{t_2}\int_{\mathbb T^2} H_{\mu_t}\cdot {\bf v}\zeta \,d\lambda_tdt.
  \end{split}
\end{equation}
\end{proposition}
\begin{proof}
Let $\rho_\delta$ be an even space–time mollifier, and write
\begin{equation}\label{liu-revision-equ5}
    \mathbf v^\delta=\rho_\delta*\mathbf v,
\qquad
\eta_\delta=\rho_\delta*(\zeta\mathbf v^\delta).
\end{equation}
First take $\zeta$ smooth and compactly supported in $\mathbb T^2\times(0,T)$, so that $\eta_\delta$ is an admissible test for sufficiently small $\delta$.
The crucial calculation uses only symmetry of convolution:

\[
\begin{aligned}
-\iint\mathbf v\cdot\partial_t\eta_\delta\,dx\,dt
&=-\iint\mathbf v^\delta\cdot
\partial_t(\zeta\mathbf v^\delta)\,dx\,dt\\
&=-\iint|\mathbf v^\delta|^2\partial_t\zeta\,dx\,dt
-\frac12\iint\zeta\,\partial_t|\mathbf v^\delta|^2\,dx\,dt\\
&=-\frac12\iint|\mathbf v^\delta|^2\partial_t\zeta\,dx\,dt\\
&\xrightarrow{\delta\to 0}
-\frac12\iint|\mathbf v|^2\partial_t\zeta\,dx\,dt
\end{aligned}
\]
by \eqref{regularity v}.
By the same approximation argument as in \eqref{mollif-v},
$\eta_\delta\to\zeta\mathbf v$ strongly in $L^2(d\lambda_tdt)$.
Since $H_{\mu_t}\in L^2(d\lambda_tdt)$, the Cauchy--Schwarz inequality  implies 
\[
\left|\iint H_{\mu_t}\cdot
(\eta_\delta-\zeta\mathbf v)\,d\lambda_t\,dt\right|
\le
\|H_{\mu_t}\|_{L^2(d\lambda_tdt)}
\|\eta_\delta-\zeta\mathbf v\|_{L^2(d\lambda_tdt)}
\xrightarrow{\delta\to 0}0.
\]
So we are able to pass to the limit
\[
\iint H_{\mu_t}\cdot\eta_\delta\,d\lambda_t\,dt
\xrightarrow{\delta\to 0}
\iint\zeta H_{\mu_t}\cdot\mathbf v\,d\lambda_t\,dt.
\]
By \eqref{regularity v} and  Ladyzhenskaya inequalities in 2d, we have $\mathbf v\in L^4_{x,t}$, and therefore
\[\iint(-\nabla\mathbf v+\mathbf v\otimes\mathbf v):\nabla\eta_\delta\,
dx\,dt\xrightarrow{\delta\to 0}\iint(-\nabla\mathbf v+\mathbf v\otimes\mathbf v):\nabla(\zeta{\bf v})\,
dx\,dt.\]
 Concerning the term involving pressure,
 \[
\begin{aligned}
\iint\Pi \,{\rm div}\,\eta_\delta \,dx\,dt
&=\iint(\Pi* \rho_\delta)\, {\rm div} (\zeta \mathbf v^\delta)\,dx\,dt\\
&=\iint (\Pi* \rho_\delta)\, \nabla\zeta\cdot \mathbf v^\delta\xrightarrow{\delta\to 0} \iint \Pi\, \nabla\zeta\cdot\mathbf v.
\end{aligned}
\]
Using the above convergence, we can  choose  $\eta_\delta$ as test function in 
\eqref{weak-NS} and calculate the $\delta$-limit by integration by parts and using ${\rm div}\,{\bf v}=0$. Suitable temporal cutoffs approximating $\mathbf1_{[t_1,t_2]}(t)$ give the
desired equality \eqref{weNS}.
\end{proof}

\begin{proof}
[Proof of Brakke's inequality  \eqref{brakke-inequ}]
By applying \eqref{localized-coupled-energy} to smooth temporal
cutoffs multiplied by a
nonnegative test function $\zeta\in C_c^\infty(\mathbb T^2\times [0,\infty))$,
we obtain the following hybrid  localized energy  inequality for a.e. $0<t_1<t_2<\infty$:
\begin{equation}\label{localized-coupled-energy-interval}
\begin{aligned}
\Big(
&\frac12\int_{\mathbb T^2}\zeta|\mathbf v|^2\,dx
+\int_{\mathbb T^2}\zeta\,d\lambda_t
\Big)\Big|_{t=t_1}^{t_2}\\
&+\int_{t_1}^{t_2}\int_{\mathbb T^2}
\zeta|\nabla\mathbf v|^2\,dx\,dt
+\int_{t_1}^{t_2}\int_{\mathbb T^2}
\zeta|H_{\mu_t}|^2\,d\lambda_t\,dt
\\
&\le
\int_{t_1}^{t_2}\int_{\mathbb T^2}
\left[\frac12|\mathbf v|^2(\partial_t\zeta+\Delta\zeta)
+\big(\frac12|\mathbf v|^2+\Pi\big)
\mathbf v\cdot\nabla\zeta\right]
\,dx\,dt
\\
&\quad\quad
+
\int_{t_1}^{t_2}\int_{\mathbb T^2}
  \partial_t \zeta +(H_{\mu_t}+P^\perp\mathbf v)\cdot\nabla\zeta
\,d\lambda_t\,dt.
\end{aligned}
\end{equation}
Substituting \eqref{weNS} into \eqref{localized-coupled-energy-interval}, we have  
\begin{equation} 
    \left.\int_{\mathbb T^2}\zeta\,d\mu_t\right|_{t=t_1}^{t_2}\leq \int_{t_1}^{t_2}\int_{\mathbb T^2}\partial_t\zeta+(H_{\mu_t}+P^\perp{\bf v})\cdot(\nabla\zeta-H_{\mu_t}\zeta)\,d\mu_t dt,
\end{equation}
where we used $\lambda_t=\sigma^{-1}\mu_t$ and  $H_{\mu_t}\cdot{\bf v}=H_{\mu_t}\cdot P^\perp{\bf v}$ due
to Brakke's perpendicularity theorem \cite[Chapter 5]{Brakke}.  
This proves \eqref{brakke-inequ}.
\end{proof}


\section{Proof of  BV identity \eqref{def:bv-flow}}
\label{sec:BV-flow}\label{section-bv}

We now show that the limiting phase satisfies the BV motion law associated
with the generalized velocity $\mathcal V$.  The argument is similar to the
proof of Tashiro \cite{Tashiro,Tashiro26}. Since the
present setting differs slightly, we give a self-contained proof for the convenience of the reader.

Let  $E$ be the measurable set obtained in Proposition
\ref{MR-compactness}, and let $E(t)=\{x\in\mathbb T^2:(x,t)\in E\}$ be its $t$-slice.  As before,
$
d\mu=d\mu_t\,dt,
$
and we write $\nabla'=(\nabla,\partial_t)$ for the space--time gradient.

The first proposition rules out jumps of the phase in the time direction.

\begin{proposition}[Space--time absolute continuity]\label{spacetime-AC}
The set $E$ has locally finite perimeter in
$\mathbb T^2\times(0,T)$ and
\begin{equation}\label{spacetime-AC-eq}
\|\nabla'\mathcal X_{E}\|\ll\mu.
\end{equation}
\end{proposition}

\begin{proof}
With 
$
\Phi(s):=\int_{-1}^s\sqrt{2W(r)}\,dr,$ we have 
$
\Phi(-1)=0$ and  $\Phi(1)=\sigma.
$
By \eqref{uk-conv-char1} 
and $|u_{\epsilon_k}|\le1$, we deduce  that
$$
\Phi(u_{\epsilon_k})
\to
\sigma\mathcal X_{E}
\qquad\text{strongly in }L^1(\mathbb T^2\times(0,T)).
$$

On $\{|\nabla u_{\epsilon_k}|>0\}$ set
$$
f_k
:=
\left(
1+
\frac{|\partial_tu_{\epsilon_k}|^2}
{|\nabla u_{\epsilon_k}|^2}
\right)^{1/2},
$$
and set $f_k=0$ on $\{|\nabla u_{\epsilon_k}|=0\}$.  By
\eqref{def: mu-tilde} and \eqref{MR-uniform-bounds}, we find 
$$
\begin{aligned}
\int f_k^2\,d\widetilde\mu_{\epsilon_k}
&=
\int_{\{|\nabla u_{\epsilon_k}|>0\}}
\left(
\epsilon_k|\nabla u_{\epsilon_k}|^2
+
\epsilon_k|\partial_tu_{\epsilon_k}|^2
\right)dx\,dt
 \le C(T).
\end{aligned}
$$  
Since
$\widetilde\mu_{\epsilon_k}\to\mu$ by \eqref{tilde-mu-limit}, Lemma
\ref{thmB.3} gives, after passing to a further subsequence which is not
relabeled, a nonnegative $f\in L^2(\mu)$ such that
$$
(\widetilde\mu_{\epsilon_k},f_k)\to(\mu,f)
$$
in the sense of measure--function pairs.

On $\{|\nabla u_{\epsilon_k}|>0\}$, set
$$
A_k:=\sqrt{\epsilon_k}|\nabla u_{\epsilon_k}|,
\qquad
B_k:=\sqrt{\frac{2W(u_{\epsilon_k})}{\epsilon_k}}.
$$
Then
$$
|\nabla'\Phi(u_{\epsilon_k})|=f_kA_kB_k,
\qquad
d\widetilde\mu_{\epsilon_k}=A_k^2\,dx\,dt,
$$
and
$$
|A_k-B_k|^2
\le
|A_k^2-B_k^2|
=
2\left|
\frac{\epsilon_k}{2}|\nabla u_{\epsilon_k}|^2
-\frac{W(u_{\epsilon_k})}{\epsilon_k}
\right|.
$$
Thus, for every nonnegative
$\phi\in C_c^0(\mathbb T^2\times(0,T))$,
$$
\begin{aligned}
&\left|
\int_{\{|\nabla u_{\epsilon_k}|>0\}}
\phi f_k(A_kB_k-A_k^2)\,dx\,dt
\right|
\\
&\qquad\le
\left(
\int\phi f_k^2\,d\widetilde\mu_{\epsilon_k}
\right)^{1/2}
\left(
2\|\phi\|_{C^0}
|\xi_{\epsilon_k}|(\mathbb T^2\times(0,T))
\right)^{1/2}
\to0
\end{aligned}
$$
by \eqref{discrepancy-vanish}.
On $\{|\nabla u_{\epsilon_k}|=0\}$,
$$
|\nabla'\Phi(u_{\epsilon_k})|
=
|\partial_tu_{\epsilon_k}|
\sqrt{2W(u_{\epsilon_k})}.
$$
Hence
$$
\begin{aligned}
&\int_{\{|\nabla u_{\epsilon_k}|=0\}}
\phi|\nabla'\Phi(u_{\epsilon_k})|\,dx\,dt
\\
&\qquad\le
\|\phi\|_{C^0}
\left(
\int\epsilon_k|\partial_tu_{\epsilon_k}|^2\,dx\,dt
\right)^{1/2}
\left(
2|\xi_{\epsilon_k}|(\mathbb T^2\times(0,T))
\right)^{1/2}
\to0,
\end{aligned}
$$
because on this set
$W(u_{\epsilon_k})/\epsilon_k=|\xi_{\epsilon_k}|$.
We have therefore proved
$$
\left|
\int\phi|\nabla'\Phi(u_{\epsilon_k})|\,dx\,dt
-
\int\phi f_k\,d\widetilde\mu_{\epsilon_k}
\right|
\to0.
$$

By lower semicontinuity of the total variation and the
measure--function convergence,
$$
\begin{aligned}
\sigma\int\phi\,d\|\nabla'\mathcal X_{E}\|
&\le
\liminf_{k\to\infty}
\int\phi|\nabla'\Phi(u_{\epsilon_k})|\,dx\,dt
 =
\int\phi f\,d\mu.
\end{aligned}
$$
Since this holds for every nonnegative $\phi\in C_c^0$, we obtain
\eqref{spacetime-AC-eq}.
\end{proof}

The preceding argument is the Allen--Cahn mechanism used in
\cite[Proposition 3.4]{Tashiro}; compare also
\cite[Proposition 8.4]{MugnaiRoger}.  In particular, it does not use a
strong smooth approximation of the fluid velocity.

We also need an upper density bound for the space--time measure.

\begin{lemma}[Space--time upper density]\label{spacetime-density}
For $\mu$-a.e.\ $(x,t)\in\mathbb T^2\times(0,T)$,
$$
\limsup_{r\downarrow0}
\frac{\mu(B_r^3(x,t))}{r^2}<\infty.
$$
Consequently,
$$
\mu\ll\mathcal H^2
\qquad\text{on }\mathbb T^2\times(0,T).
$$
\end{lemma}

\begin{proof}
By \eqref{limit-density-L2},
$$
D(t)
=
\sup_{x\in\mathbb T^2,\ 0<r\leq 1}
\frac{\mu_t(B_r(x))}{r}
\in L^2(0,T).
$$
If $t$ is a Lebesgue point of $D^2$ and $r>0$ is sufficiently small, then
$$
\begin{aligned}
\mu(B_r^3(x,t))
&\le
\int_{t-r}^{t+r}\mu_s(B_r(x))\,ds
 \le
r\int_{t-r}^{t+r}D(s)\,ds.
\end{aligned}
$$
Therefore
$$
\frac{\mu(B_r^3(x,t))}{r^2}
\le
\left(
\frac{2}{r}\int_{t-r}^{t+r}D(s)^2\,ds
\right)^{1/2},
$$
whose limit superior is finite.  Since the exceptional set of times has
Lebesgue measure zero and
$\operatorname*{ess\,sup}_{0<t<T}\mu_t(\mathbb T^2)<\infty$, the upper
two-dimensional density is finite for $\mu$-a.e.\ $(x,t)$.

To obtain the last assertion, for $j\in\mathbb N$ let
$$
G_j
=
\left\{
z:
\limsup_{r\downarrow0}\frac{\mu(B_r^3(z))}{r^2}\le j
\right\}.
$$
The preceding argument gives
$\mu((\mathbb T^2\times(0,T))\setminus\bigcup_jG_j)=0$.
The standard density comparison theorem gives
$\mu(A\cap G_j)\le Cj\mathcal H^2(A)$ for every Borel set $A$.
Thus $\mathcal H^2(A)=0$ implies $\mu(A)=0$.   \end{proof}

We next combine Proposition \ref{spacetime-AC}, Lemma
\ref{spacetime-density}, and the $L^2$-flow property.

\begin{proposition}\label{spacetime-tangent}
For $\mathcal H^2$-a.e.\
$(x,t)\in\partial^*E$, the following properties hold:
\begin{enumerate}
\item
the tangent plane $T_{(x,t)}\mu$ exists and
$$
T_{(x,t)}\mu
=
T_{(x,t)}(\partial^*E),
\qquad
\binom{\mathcal V(x,t)}{1}
\in T_{(x,t)}\mu;
$$

\item
$x\in\partial^*E(t)$ and
$$
T_x\mu_t=T_x(\partial^*E(t));
$$

\item
if $\nu_{E}$ is the outer unit normal to $E$ and
$\mathbf p:\mathbb R^2\times\mathbb R\to\mathbb R^2$ denotes the spatial
projection, then
$$
\mathbf p(\nu_{E}(x,t))\ne0,
\qquad
\nu_{E(t)}(x)
=
\frac{\mathbf p(\nu_{E}(x,t))}
{|\mathbf p(\nu_{E}(x,t))|};
$$

\item
$$
T_x(\partial^*E(t))\times\{0\}
\subset
T_{(x,t)}(\partial^*E).
$$
\end{enumerate}
\end{proposition}

\begin{proof}
Set
$$
M:=\partial^*E\subset\mathbb T^2\times(0,T).
$$
Since $E$ has locally finite perimeter in space--time, De Giorgi's
structure theorem gives
\begin{equation}\label{DE-giorgi-space-time}
\|\nabla'\mathcal X_E\|
=
\mathcal H^2\llcorner M,
\end{equation}
and $M$ is countably $2$-rectifiable.

We divide the proof into several steps.

\medskip
\noindent
\emph{Step 1: the tangent plane of the space--time measure on $M$.}

By Proposition \ref{spacetime-AC} and
\eqref{DE-giorgi-space-time},
$$
\mathcal H^2\llcorner M\ll\mu.
$$
On the other hand, Lemma \ref{spacetime-density} gives
$$
\mu\ll\mathcal H^2.
$$
In particular, after restriction to $M$,
$$
\mu\llcorner M
\ll
\mathcal H^2\llcorner M
\ll
\mu\llcorner M.
$$
Since $\mathcal H^2\llcorner M$ is locally finite, the
Radon--Nikod\'ym theorem yields a function
$$
\theta\in L^1_{\rm loc}
(\mathcal H^2\llcorner M)
$$
such that
\begin{equation}\label{mu-on-M}
\mu\llcorner M
=
\theta\,\mathcal H^2\llcorner M,
\end{equation}
where
\begin{equation}\label{theta-positive}
0<\theta(x,t)<\infty
\qquad
\text{for $\mathcal H^2$-a.e.\ $(x,t)\in M$}.
\end{equation}

Since $M$ is countably $2$-rectifiable, for
$\mathcal H^2$-a.e.\ $z=(x,t)\in M$ the approximate tangent plane
$T_zM$ exists.  At such a point which is also a Lebesgue point of
$\theta$, one has
\begin{equation}\label{restriction-tangent}
r^{-2}
(\eta_{z,r})_\#
(\mu\llcorner M)
\longrightarrow
\theta(z)\,
\mathcal H^2\llcorner T_zM
\end{equation}
as Radon measures on $\mathbb R^3$, where
$$
\eta_{z,r}(y,s)
=
\left(
\frac{y-x}{r},\frac{s-t}{r}
\right).
$$

It remains to show that the part of $\mu$ outside $M$ does not
contribute to this blow-up.  Set
$$
\mu_1:=\mu\llcorner M,
\qquad
\mu_2:=\mu\llcorner
\big((\mathbb T^2\times(0,T))\setminus M\big).
$$
The measures $\mu_1$ and $\mu_2$ are mutually singular.  By the
differentiation theorem for Radon measures,
\begin{equation}\label{differentiate-singular}
\lim_{r\downarrow0}
\frac{\mu_2(B_r^3(z))}
{\mu_1(B_r^3(z))}
=0
\end{equation}
for $\mu_1$-a.e.\ $z$.  This is the differentiation result used in
\cite[Theorem 3.5]{SimonGMT}; compare the corresponding step in
\cite[Proposition 3.5]{Tashiro} and
\cite[Lemma 4.7]{StuvardTonegawa}.

By \eqref{mu-on-M}--\eqref{theta-positive}, the exceptional set in
\eqref{differentiate-singular} is also
$\mathcal H^2$-null on $M$.  Moreover, at the points under
consideration,
$$
\lim_{r\downarrow0}
\frac{\mu_1(B_r^3(z))}{\pi r^2}
=
\theta(z).
$$
Hence \eqref{differentiate-singular} gives
\begin{equation}\label{outside-zero-density}
\lim_{r\downarrow0}
\frac{
\mu(B_r^3(z)\setminus M)
}{r^2}
=0
\qquad
\text{for $\mathcal H^2$-a.e.\ $z\in M$}.
\end{equation}
Combining \eqref{restriction-tangent} and
\eqref{outside-zero-density}, we obtain
$$
r^{-2}(\eta_{z,r})_\#\mu
\longrightarrow
\theta(z)\,
\mathcal H^2\llcorner T_zM.
$$
Thus the tangent plane $T_z\mu$ exists and
\begin{equation}\label{space-time-tangent-identification}
T_z\mu=T_zM
=
T_{(x,t)}(\partial^*E)
\end{equation}
for $\mathcal H^2$-a.e.\ $z=(x,t)\in M$.

\medskip
\noindent
\emph{Step 2: the space--time tangent contains the generalized
velocity.}

Set
$$
d\widehat\mu_t:=\sigma^{-1}d\mu_t,
\qquad
d\widehat\mu:=d\widehat\mu_t\,dt
=\sigma^{-1}d\mu.
$$
By Proposition \ref{MR-compactness},
$\{\widehat\mu_t\}_{t\in(0,T)}$ is an $L^2$-flow with generalized
velocity $\mathcal V$.

We use \cite[Proposition 3.3]{MugnaiRoger}.  In the notation of
Mugnai--R\"oger, where space--time points are ordered as $(t,x)$,
that proposition states that
$$
\binom{1}{\mathcal V(x,t)}
\in T_{(t,x)}\widehat\mu
$$
at $\widehat\mu$-a.e.\ point where the tangent plane exists.
With our convention $(x,t)$ for space--time coordinates, this becomes
\begin{equation}\label{MR-tangent-property}
\binom{\mathcal V(x,t)}{1}
\in T_{(x,t)}\widehat\mu
=
T_{(x,t)}\mu.
\end{equation}
Multiplication of a measure by the positive constant $\sigma^{-1}$
does not change its tangent plane.

Since
$\mathcal H^2\llcorner M\ll\mu,$
a $\mu$-null set is $\mathcal H^2$-null on $M$.  Thus, combining
\eqref{space-time-tangent-identification} and
\eqref{MR-tangent-property},
\begin{equation}\label{velocity-tangent-M}
\binom{\mathcal V(x,t)}{1}
\in T_{(x,t)}M
\end{equation}
for $\mathcal H^2$-a.e.\ $(x,t)\in M$.

Let
$
q:\mathbb T^2\times(0,T)\to(0,T),
q(x,t)=t.
$
At a point where $T_{(x,t)}M$ exists, write
$$
\nabla^M q
=
P_{T_{(x,t)}M}(0,0,1).
$$
From \eqref{velocity-tangent-M},
$
\begin{aligned}
1
&=
(0,0,1)\cdot
\binom{\mathcal V}{1}
=
\nabla^M q\cdot
\binom{\mathcal V}{1}.
\end{aligned}
$
Therefore
\begin{equation}\label{q-positive}
|\nabla^M q|
\ge
\frac1{\sqrt{1+|\mathcal V|^2}}
>0
\end{equation}
for $\mathcal H^2$-a.e.\ point of $M$.

The strict positivity in \eqref{q-positive} is important below.

\medskip
\noindent
\emph{Step 3: slicing the reduced boundary.}

We use the slicing theorem for sets of finite perimeter,
\cite[Theorem 18.11]{Maggi}.  Applied in local coordinates on
$\mathbb T^2\times(0,T)$ to the time projection $q$, it gives a set
$Z_1\subset(0,T)$ with $\mathcal L^1(Z_1)=0$ such that, for every
$t\notin Z_1$,
$E(t)$ has locally finite perimeter in $\mathbb T^2$ and
\begin{equation}\label{slicing-boundaries}
\mathcal H^1
\left(
M_t\mathbin{\triangle}\partial^*E(t)
\right)
=0,
\qquad
M_t:=\{x:(x,t)\in M\},
\end{equation}
and, for $\mathcal H^1$-a.e.\ $x\in M_t$,
\begin{equation}\label{slicing-normal}
\mathbf p(\nu_E(x,t))\ne0,
\qquad
\nu_{E(t)}(x)
=
\frac{\mathbf p(\nu_E(x,t))}
{|\mathbf p(\nu_E(x,t))|}.
\end{equation}
These are precisely the slicing conclusions recorded in
\cite[Theorem 18.11]{Maggi}; they are also stated explicitly in
\cite[Lemma 3.8]{Tashiro}.

We now explain carefully why these statements, which initially hold
for a.e.\ $t$ and $\mathcal H^1$-a.e.\ point of the slice, hold for
$\mathcal H^2$-a.e.\ point of $M$.

Let $B_1\subset M$ be the set of points $(x,t)$ for which either
$t\in Z_1$, or $x\notin\partial^*E(t)$, or
\eqref{slicing-normal} fails.  For every $t\notin Z_1$,
\eqref{slicing-boundaries}--\eqref{slicing-normal} give
$\mathcal H^1((B_1)_t)=0.$
The coarea formula on the countably $2$-rectifiable set $M$ gives
$$
\begin{aligned}
\int_{B_1}
|\nabla^M q|\,d\mathcal H^2
&=
\int_0^T
\mathcal H^1((B_1)_t)\,dt
=0.
\end{aligned}
$$
Since $|\nabla^M q|>0$ for $\mathcal H^2$-a.e.\ point of $M$ by
\eqref{q-positive}, it follows that
$
\mathcal H^2(B_1)=0.
$
Therefore, for $\mathcal H^2$-a.e.\ $(x,t)\in M$,
\begin{equation}\label{slice-conclusions-spacetime}
x\in\partial^*E(t),
\qquad
\mathbf p(\nu_E(x,t))\ne0,
\qquad
\nu_{E(t)}(x)
=
\frac{\mathbf p(\nu_E(x,t))}
{|\mathbf p(\nu_E(x,t))|}.
\end{equation}
This proves the first assertion in part (2) and all of part (3).

\medskip
\noindent
\emph{Step 4: equality of the spatial tangent lines.}

By Proposition \ref{MR-compactness}, there exists a set
$Z_2\subset(0,T)$ of Lebesgue measure zero such that, for
$t\notin Z_2$,
$\sigma^{-1}\mu_t$ is an integral one-dimensional varifold and
\begin{equation}\label{phase-domination-current}
\mathcal H^1\llcorner\partial^*E(t)
\le
\sigma^{-1}\mu_t.
\end{equation}

Fix such a time $t$.  Since $\sigma^{-1}\mu_t$ is integral, there
exist a countably $1$-rectifiable set $N_t\subset\mathbb T^2$ and an
integer-valued function
$\vartheta_t\ge1$ such that
$$
\sigma^{-1}\mu_t
=
\vartheta_t\,
\mathcal H^1\llcorner N_t.
$$
From \eqref{phase-domination-current},
$$
\mathcal H^1
\big(
\partial^*E(t)\setminus N_t
\big)
=0.
$$
Both $N_t$ and $\partial^*E(t)$ are countably $1$-rectifiable.
At $\mathcal H^1$-a.e.\ point of their intersection their
approximate tangent lines therefore agree.  Indeed, this follows
from uniqueness of the approximate tangent line for the rectifiable
set
$
N_t\cup\partial^*E(t).
$
Consequently,
\begin{equation}\label{slice-tangent-equality}
T_x\mu_t
=
T_xN_t
=
T_x(\partial^*E(t))
\end{equation}
for $\mathcal H^1$-a.e.\ $x\in\partial^*E(t)$.

Again we have to transfer this fixed-time conclusion to a
space--time statement.  Let $B_2\subset M$ be the set of points
$(x,t)$ for which \eqref{slice-tangent-equality} fails, including
the exceptional times in $Z_1\cup Z_2$.  By
\eqref{slicing-boundaries} and \eqref{slice-tangent-equality},
$
\mathcal H^1((B_2)_t)=0
$
for a.e.\ $t$.  Hence the coarea formula gives
$$
\int_{B_2}
|\nabla^M q|\,d\mathcal H^2
=
\int_0^T
\mathcal H^1((B_2)_t)\,dt
=0.
$$
Using \eqref{q-positive} once more,
$
\mathcal H^2(B_2)=0.
$
Thus
\begin{equation}\label{spacetime-to-slice-tangent}
T_x\mu_t
=
T_x(\partial^*E(t))
\end{equation}
for $\mathcal H^2$-a.e.\ $(x,t)\in M$.
This completes part (2).

\medskip
\noindent
\emph{Step 5: inclusion of the spatial tangent line in the
space--time tangent plane.}

Take a point $(x,t)\in M$ at which all the conclusions proved above
hold, and let
$$
z\in T_x(\partial^*E(t)).
$$
By definition of the reduced boundary,
$$
z\cdot\nu_{E(t)}(x)=0.
$$
By \eqref{slice-conclusions-spacetime},
$$
\mathbf p(\nu_E(x,t))
=
|\mathbf p(\nu_E(x,t))|
\,\nu_{E(t)}(x).
$$
Therefore
$$
\begin{aligned}
(z,0)\cdot\nu_E(x,t)
&=
z\cdot\mathbf p(\nu_E(x,t))
\\
&=
|\mathbf p(\nu_E(x,t))|
\,z\cdot\nu_{E(t)}(x)
\\
&=0.
\end{aligned}
$$
Since
$$
T_{(x,t)}(\partial^*E)
=
\nu_E(x,t)^\perp,
$$
we conclude that
$$
(z,0)
\in
T_{(x,t)}(\partial^*E).
$$
As this holds for every
$z\in T_x(\partial^*E(t))$,
$$
T_x(\partial^*E(t))\times\{0\}
\subset
T_{(x,t)}(\partial^*E).
$$
This proves part (4), and the proof is complete.
\end{proof}

This proposition is the analogue of
\cite[Proposition 3.5]{Tashiro}; see also
\cite[Lemma 4.7]{StuvardTonegawa}.  The only flow property used in its
geometric part is the $L^2$-flow tangent relation above.

\begin{proposition}[BV motion law]\label{BV-motion}
There is a representative of $t\mapsto\mathcal X_{E(t)}$ which belongs to
$$
C^{0,1/2}
\big([0,T];L^1(\mathbb T^2)\big)
$$
and such that, for every $0\leq t_1<t_2\leq T$ and every
$\phi\in C^1(\mathbb T^2\times[0,T])$,
\begin{equation}\label{BV-total-velocity}
\begin{aligned}
&\left.
\int_{E(t)}\phi(x,t)\,dx
\right|_{t=t_1}^{t_2}
\\
&\qquad=
\int_{t_1}^{t_2}\int_{E(t)}
\partial_t\phi\,dx\,dt
+
\int_{t_1}^{t_2}\int_{\partial^*E(t)}
\phi\,\mathcal V\cdot\nu_{E(t)}
\,d\mathcal H^1\,dt.
\end{aligned}
\end{equation}
Moreover,
$$
\int_0^T\int_{\partial^*E(t)}
|\mathcal V\cdot\nu_{E(t)}|^2
\,d\mathcal H^1\,dt<\infty.
$$
\end{proposition}

\begin{proof}
Fix a point $(x,t)\in\partial^*E$ at which all conclusions of
Proposition \ref{spacetime-tangent} hold.  Since
$$
\mathcal V(x,t)\perp T_x\mu_t
=
T_x(\partial^*E(t)),
$$
the two-dimensional tangent plane is
$$
T_{(x,t)}(\partial^*E)
=
\big(T_x(\partial^*E(t))\times\{0\}\big)
\oplus
\operatorname{span}
\left\{
\binom{\mathcal V(x,t)}{1}
\right\}.
$$
The orientation in Proposition \ref{spacetime-tangent} then gives
\begin{equation}\label{space-time-normal}
\nu_{E}(x,t)
=
\frac{1}{\sqrt{1+|\mathcal V(x,t)|^2}}
\binom{
\nu_{E(t)}(x)
}{
-\mathcal V(x,t)\cdot\nu_{E(t)}(x)
}.
\end{equation}
In particular, for the time projection $q(x,t)=t$,
$$
|\nabla^{\partial^*E}q|
=
\frac{1}{\sqrt{1+|\mathcal V|^2}}.
$$

Let $\phi\in C_c^1(\mathbb T^2\times(0,T))$.  The Gauss--Green formula
in space--time and \eqref{space-time-normal} give
$$
\begin{aligned}
\int_{\mathbb T^2\times(0,T)}
\mathcal X_{E}\partial_t\phi\,dx\,dt
&=
\int_{\partial^*E}
\phi\,(\nu_{E})_t\,d\mathcal H^2
\\
&=
-\int_{\partial^*E}
\phi\,
(\mathcal V\cdot\nu_{E(t)})
|\nabla^{\partial^*E}q|
\,d\mathcal H^2.
\end{aligned}
$$
Applying the coarea formula on $\partial^*E$ yields
\begin{equation}\label{BV-distribution}
\int_{\mathbb T^2\times(0,T)}
\mathcal X_{E}\partial_t\phi\,dx\,dt
=
-\int_0^T\int_{\partial^*E(t)}
\phi\,\mathcal V\cdot\nu_{E(t)}
\,d\mathcal H^1\,dt.
\end{equation}
Here the exceptional subsets in Proposition \ref{spacetime-tangent} have
zero contribution by the coarea formula.

Equation \eqref{BV-distribution} first gives
\eqref{BV-total-velocity} for a.e.\ pair $t_1<t_2$.  Furthermore,
by \eqref{phase-measure-domination},
$$
\begin{aligned}
&\int_0^T\int_{\partial^*E(t)}
|\mathcal V\cdot\nu_{E(t)}|^2\,d\mathcal H^1\,dt
\\
&\qquad\le
\sigma^{-1}
\int_0^T\int_{\mathbb T^2}
|\mathcal V|^2\,d\mu_t\,dt
<\infty.
\end{aligned}
$$
Also,
$$
\operatorname*{ess\,sup}_{0<t<T}
\mathcal H^1(\partial^*E(t))
\le
\sigma^{-1}
\operatorname*{ess\,sup}_{0<t<T}\mu_t(\mathbb T^2)
<\infty.
$$
Hence, for a.e.\ $t_1<t_2$,
$$
\begin{aligned}
\|\mathcal X_{E(t_2)}-\mathcal X_{E(t_1)}\|_{L^1}
&\le
\int_{t_1}^{t_2}\int_{\partial^*E(t)}
|\mathcal V|\,d\mathcal H^1\,dt
\\
&\le C(T)|t_2-t_1|^{1/2}.
\end{aligned}
$$
Thus $t\mapsto\mathcal X_{E(t)}$ has a unique
$L^1$-continuous and $1/2$-H\"older representative
on $[0,T]$.  With this
representative, \eqref{BV-total-velocity} extends by continuity to every
$0\leq t_1<t_2\leq T$.
\end{proof}

  \section{A remark on the Intrinsic mean curvature}
 
The following is useful when the limiting energy measure has a higher
multiplicity.  It is special to the present two-dimensional setting because
R\"oger's locality theorem \cite[Proposition 3.1]{RogerStefan} requires an exponent strictly larger than the
varifold dimension. 

\begin{proposition}[R\"oger's locality of genralized mean curvature]\label{intrinsic-curvature}
For a.e.\ $t>0$ there exists an $\mathcal H^1$-a.e.\ defined vector field
$H_{E(t)}$ on $\partial^*E(t)$, depending only on the reduced phase boundary,
such that
\[
H_{\mu_t}=H_{E(t)}
\qquad\mathcal H^1\text{-a.e.\ on }\partial^*E(t).
\]
More generally, if $\lambda$ is any integral one-dimensional varifold with
locally bounded first variation satisfying
$$
\mathcal H^1(\partial^*E(t)\setminus\operatorname{spt}\lambda)=0,
\qquad
H_\lambda\in L^2_{\mathrm{loc}}(\lambda;\mathbb R^2),
$$
then
$$
H_\lambda=H_{E(t)}
\qquad\mathcal H^1\text{-a.e.\ on }\partial^*E(t).
$$
In particular, the generalized mean curvature on the reduced phase boundary is
independent of the multiplicity and of the choice of an ambient integral
$L^2$-curvature varifold containing that boundary.
\end{proposition}

\begin{proof}
Fix a time for which Proposition \ref{MR-compactness} holds.  By
\eqref{phase-varifold},
$$
\mathcal H^1\llcorner\partial^*E(t)\le\sigma^{-1}\mu_t,
$$
so $\partial^*E(t)$ is contained in $\operatorname{spt}\mu_t$ up to an
$\mathcal H^1$-null set.  The varifold $\sigma^{-1}\mu_t$ is integral and has
generalized mean curvature $H_{\mu_t}\in L^2_{\mathrm{loc}}(\mu_t)$.  Since the
ambient dimension is two,   the hypothesis in \cite[Proposition 3.1]{RogerStefan} is fulfilled.
That proposition states that if two integral $(n-1)$-varifolds contain the
same reduced boundary and have generalized mean curvature in $L^s_{\mathrm{loc}}$
with $s>n-1$ and $s\ge2$, then their generalized mean curvatures coincide
$\mathcal H^{n-1}$-a.e.\ on that reduced boundary.  Applying it to
$\sigma^{-1}\mu_t$ and any admissible $\lambda$ proves the asserted
independence.

The vector field $H_{E(t)}$ is therefore defined to be the common value.  The
proof of \cite[Lemma 3.2]{RogerStefan} gives a more intrinsic description: on
each measurable graph $y\mapsto(y,\psi(y))$ covering the reduced boundary,
$\psi$ is twice approximately differentiable at a.e.\ relevant point and
$$
H_{E(t)}(y,\psi(y))
=
\left(\frac{\psi'}{\sqrt{1+|\psi'|^2}}\right)'
\frac{(-\psi',1)}{\sqrt{1+|\psi'|^2}},
$$
where the derivatives are approximate derivatives.  Hence the common
curvature is determined by the phase boundary itself.
\end{proof}

\begin{remark}
On a space-time region $U\subset\mathbb T^2\times[0,\infty)$, if the limiting measure $\sigma^{-1}\mu_t$ is of unit density for a.e.~$t$, then
$$
\sigma^{-1}\mu_t
=\mathcal H^1\llcorner\partial^*E(t).
$$
This can be proved by examining the proof of integrality 
in \cite{MugnaiRoger} or \cite{HT}. 
Thus, in this case in $U$, \eqref{velocity-law} becomes on the reduced phase boundary
$$
\mathcal V=H_{E(t)}+P^\perp\mathbf v.
$$
Thus both the curvature and the motion law can be stated
entirely in terms of the phase boundary and the fluid velocity in $U$.
\end{remark}

 {\bf Knowledgement} 
 Yuning Liu  acknowledges the support of NSF of China under Grant 12571224.
Yoshihiro Tonegawa acknowledges the support of the JSPS Grant-in-aid for scientific research 23H00085.

\medskip
\noindent

{\bf AI assistance statement}. During the preparation of this manuscript, the authors used OpenAI models for exploratory discussions, mathematical checks, and manuscript editing. The authors independently verified the mathematical content and take full responsibility for the final manuscript.

\end{document}